\documentclass{amsart}

\usepackage{amsmath, amssymb}
\usepackage{amsthm}
\usepackage[arrow,matrix,curve,cmtip,ps]{xy}

\numberwithin{equation}{section}

\allowdisplaybreaks

\newtheorem{theorem}{Theorem}[section]
\newtheorem{lemma}[theorem]{Lemma}

\newtheorem{corollary}[theorem]{Corollary}

\newtheorem*{theorem*}{Theorem}
\theoremstyle{remark}
\newtheorem{remark}[theorem]{Remark}
\newtheorem{definition}[theorem]{Definition}
\newtheorem{example}[theorem]{Example}

\theoremstyle{definition}
\newtheorem*{convention*}{Conventions}
\newcommand{\Z}{\mathbb{Z}}
\newcommand{\N}{\mathbb{N}}

\newcommand{\T}{\mathbb{T}}
\newcommand{\im}{\operatorname{im}}

\newcommand{\aut}{\operatorname{Aut}}

\newcommand{\TX}{\mathcal{T}_X}
\newcommand{\OX}{\mathcal{O}_X}
\newcommand{\OKX}{\mathcal{O}(K,X)}

\newcommand{\clspan}{\operatorname{\overline{\mathrm{span}}}}

\newcommand{\Li}{\mathcal{L}}
\newcommand{\K}{\mathcal{K}}

\newcommand{\OY}{\mathcal{O}_Y}

\newcommand{\Hi}{\mathcal{H}}
\newcommand{\B}{\mathcal{B}}

\newcommand{\Q}{\mathcal{Q}}

\begin{document}
\title[Adding Tails to $C^*$-correspondences]{Adding Tails to
$\boldsymbol{C^*}$-correspondences}
\author{Paul S. Muhly}
\author{Mark Tomforde }
\address{Department of Mathematics\\
University of Iowa\\
Iowa City\\
IA 52242-1419\\
USA}
\email{pmuhly@math.uiowa.edu}
\email{tomforde@math.uiowa.edu}
\thanks{The first author was supported by NSF Grant DMS-0070405 and the
second author was supported by NSF Postdoctoral Fellowship DMS-0201960.}
\date{\today}


\begin{abstract}
We describe a method of adding tails to $C^*$-correspondences which generalizes
the process used in the study of graph $C^*$-algebras.  We show how this
technique can be used to extend results for augmented Cuntz-Pimsner algebras to
$C^*$-algebras associated to general $C^*$-correspondences, and as an
application we prove a gauge-invariant uniqueness theorem for these algebras. 
We also define a notion of relative graph $C^*$-algebras and show that
properties of these $C^*$-algebras can provide insight and motivation for
results about relative Cuntz-Pimsner algebras.
\end{abstract}

\maketitle

%%%%%%%%%%%%%%%%%%%%%%%%%%%%%%%%%%%%%%%%%%%%%%%%%%%%%%%%%%
\section{Introduction}
%%%%%%%%%%%%%%%%%%%%%%%%%%%%%%%%%%%%%%%%%%%%%%%%%%%%%%%%%%

In \cite{Pim} Pimsner introduced a way to construct a $C^*$-algebra
$\mathcal{O}_{X}$ from a pair $(A,X)$, where $A$ is a $C^*$-algebra and $X$ is
a $C^*$-correspondence (sometimes called a Hilbert bimodule) over
$A$.  Throughout his analysis Pimsner assumed that his correspondence was full
and that the left action of $A$ on $X$ was injective.  These Cuntz-Pimsner
algebras have been found to compose a class of $C^*$-algebras that is
extraordinarily rich and includes numerous $C^*$-algebras found in the
literature: crossed products by automorphisms, crossed products by
endomorphisms, partial crossed products, Cuntz-Krieger algebras, $C^*$-algebras
of graphs with no sinks, Exel-Laca algebras, and many more. Consequently, the
study of Cuntz-Pimsner algebras has received a fair amount of attention by the
operator algebra community in recent years, and because information about
$\mathcal{O}_{X}$ is very densely codified in $(A,X)$, determining how to
extract it has been the focus of much current effort.

One interesting consequence of this effort has been the introduction of the
so-called \emph{relative Cuntz-Pimsner algebras}, denoted $\mathcal{O} (K,X)$,
that have Cuntz-Pimsner algebras as quotients. Very roughly speaking, a
relative Cuntz-Pimsner algebra arises by relaxing some of the relations that
must hold among the generators of a Cuntz-Pimsner algebra. These relations are
codified in an ideal $K$ of $A$. (The precise definition will be given
shortly.) Relative Cuntz-Pimsner algebras arise quite naturally, particularly
when trying to understand the ideal structure of a Cuntz-Pimsner algebra (See,
e.g., \cite{MS, FMR}). It turns out, in fact, that not only are Cuntz-Pimsner
algebras quotients of relative Cuntz-Pimsner algebras, but quotients of
Cuntz-Pimsner algebras are often relative Cuntz-Pimsner algebras
\cite[Theorem 3.1]{FMR}.

Although in his initial work Pimsner assumed that his $C^*$-correspondences
were full and had injective left action, in recent years there have been
efforts to remove these restrictions.  Pimsner himself described how to
deal with the case when $X$ was not full, defining the so-called
\emph{augmented Cuntz-Pimsner algebras} \cite[Remark~1.2(3)]{Pim}.  However, the
case when the left action is not injective has been more elusive.  In
\cite{FMR} it was shown that for any $C^*$-correspondence $X$ and for any ideal
$K$ of $A$ consisting of elements that act as compact operators on the left of
$X$, one may define $\mathcal{O}(K,X)$ to be a $C^*$-algebra which satisfies a
certain universal property \cite[Proposition~1.3]{FMR}.  In the case that $X$
is full with injective left action, this definition agrees with previously
defined notions of relative Cuntz-Pimsner algebras, and the
Cuntz-Pimsner algebra $\mathcal{O}_X$ is equal to $\mathcal{O}(J(X),X)$, where
$J(X)$ denotes the ideal consisting of all elements of $A$ which act on the left
of $X$ as compact operators.

In \cite{FMR} it was proposed that for a general $C^*$-correspondence $X$, the
$C^*$-algebra $\mathcal{O}(J(X),X)$ is the proper analogue of the
Cuntz-Pimsner algebra.  However, upon further analysis it seems that this is
not exactly correct.  To see why, consider the case of graph $C^*$-algebras. 
If $E = (E^0,E^1,r,s)$ is a graph, then there is a natural $C^*$-correspondence
$X(E)$ over $C_0(E^0)$ associated to $E$ (see \cite[Example~1.2]{FR}).  If $E$
has no sinks, then the $C^*$-algebra $\mathcal{O}(J(X(E)), X(E))$ is isomorphic
to the graph $C^*$-algebra $C^*(E)$.  However, when $E$ has sinks this will not
necessarily be the case.  

It is worth mentioning that graphs with sinks play an important role in the
study of graph $C^*$-algebras.  Even if one begins with a graph $E$
containing no sinks, an analysis of $C^*(E)$ will often necessitate
considering $C^*$-algebras of graphs with sinks.  For example, quotients of
$C^*(E)$ will often be isomorphic to $C^*$-algebras of graphs with sinks even
when $E$ has no sinks.  Consequently, one needs a theory that incorporates
these objects.

This deficiency in the generalization of Cuntz-Pimsner algebras was addressed by
Katsura in \cite{Kat3} and \cite{Kat4}.  If $X$ is a $C^*$-correspondence over a
$C^*$-algebra $A$ with left action $\phi : A \to \Li (X)$, then Katsura proposed
that the appropriate analogue of the Cuntz-Pimsner algebra is
$\mathcal{O}_X := \mathcal{O}(J_X,X)$, where $$J_X := \{ a \in J(X) : ab=0
\text{ for all $b \in \ker \phi$} \}.$$  (Note that when $\phi$ is injective
$J_X = J(X)$.)  It turns out that when $\phi$ is injective, $\mathcal{O}_X$
is equal to the augmented Cuntz-Pimsner algebra of $X$, and when $X$ is also
full $\mathcal{O}_X$ coincides with the Cuntz-Pimsner algebra of $X$. 
Furthermore, if $E$ is a graph (possibly containing sinks), then
$\mathcal{O}_{X(E)}$ is isomorphic to $C^*(E)$.  In addition, as with graph
algebras, the class of $\mathcal{O}_X$'s is closed under quotients by
gauge-invariant ideals.  These facts, together with the analysis described
in \cite{Kat4} and \cite{Kat5}, provide strong arguments for using
$\mathcal{O}_X := \mathcal{O}(J_X,X)$ as the analogue of the Cuntz-Pimsner
algebra.  We shall adopt this viewpoint here, and for a general
$C^*$-correspondence $X$ we define $\mathcal{O}_X := \mathcal{O}(J_X,X)$ to be
the $C^*$-algebra associated to $X$.

In this paper we shall describe a method which will allow one to ``bootstrap"
many results for augmented Cuntz-Pimsner algebras to $C^*$-algebras associated
to general correspondences.  This method is inspired by a technique from the
theory of graph $C^*$-algebras, where one can often reduce to the sinkless case
by the process of ``adding tails to sinks".  Specifically, if $E$ is a graph and
$v$ is a vertex of $E$, then by \emph{adding a tail to $v$} we mean attaching a
graph of the form
\begin{equation*}
\xymatrix{ v \ar[r] & \bullet \ar[r] & \bullet \ar[r] &
\bullet \ar[r] & \cdots\\ }
\end{equation*}
to $E$.  It is well known that if $F$ is the graph formed by
adding a tail to every sink of $E$, then $F$ is a graph with no sinks and
$C^*(E)$ is canonically isomorphic to a full corner of $C^*(F)$.  Thus in the
proofs of many theorems about graph $C^*$-algebras, one can reduce to the case
of no sinks.

In this paper we describe a generalization of this process for
$C^*$-correspondences.  More specifically, if $X$ is a
$C^*$-correspondence over a $C^*$-algebra $A$, then we describe how to
construct a $C^*$-algebra $B$ and a $C^*$-correspondence $Y$ over $B$ with the
property that the left action of $Y$ is injective and $\mathcal{O}_X$ is
canonically isomorphic to a full corner of $\mathcal{O}_Y$.  Thus many
questions about $C^*$-algebras associated to correspondences can be reduced to
questions about augmented Cuntz-Pimsner algebras, and many results
characterizing properties of augmented Cuntz-Pimsner algebras may be
easily generalized to $C^*$-algebras associated to general correspondences.  As
an application of this technique, we use it in the proof of Theorem~\ref{GIU}
to extend the Gauge-Invariant Uniqueness Theorem for augmented Cuntz-Pimsner
algebras to $C^*$-algebras of general correspondences.

This paper is organized as follows.  We begin in Section~\ref{prelim-sec} with
some preliminaries.  In Section~\ref{rel-gr-sec} we analyze graph
$C^*$-algebras in the context of Cuntz-Pimsner and relative Cuntz-Pimsner
algebras, and describe a notion of a \emph{relative graph $C^*$-algebra}.
Since graph algebras provide much of the impetus for our analysis of
$C^*$-correspondences, we examine these objects carefully in order to provide a
framework which will motivate and illuminate the results of subsequent
sections.  In Section~\ref{add-tails-sec} we describe our main result --- a
process of ``adding tails" to general $C^*$-correspondences.  We also prove
that this process preserves the Morita equivalence class of the associated
$C^*$-algebra.  In Section~\ref{GIU-sec} we provide an application of our
technique of ``adding tails" by using it to extend the Gauge-Invariant
Uniqueness Theorem for augmented Cuntz-Pimsner algebras to $C^*$-algebras
associated to general correspondences.  We also interpret this theorem in the
context of relative Cuntz-Pimsner algebras, and in Section~\ref{ideals-sec} we use it to classify the gauge-invariant ideals in $C^*$-algebras associated to certain correspondences.  Finally, we conclude in
Section~\ref{conc-sec} by discussing other possible applications of our technique.

The authors would like to thank Takeshi Katsura for pointing out an error in a
previous draft of this paper, and for many useful conversations regarding these
topics.

%%%%%%%%%%%%%%%%%%%%%%%%%%%%%%%%%%%%%%%%%%%%%%%%%%%%%%%%%%
\section{Preliminaries} \label{prelim-sec}
%%%%%%%%%%%%%%%%%%%%%%%%%%%%%%%%%%%%%%%%%%%%%%%%%%%%%%%%%%

For the most part we will use the notation and conventions of \cite{FMR},
augmenting them when necessary with the innovations of \cite{Kat3} and
\cite{Kat4}.

\begin{definition}
If $A$ is a $C^*$-algebra, then a \emph{right Hilbert $A$-module} is a Banach
space $X$ together with a right action of $A$ on $X$ and an $A$-valued inner
product $\langle \cdot , \cdot \rangle_A$ satisfying
\begin{enumerate}
\item[(i)] $\langle \xi, \eta a \rangle_A =  \langle \xi, \eta \rangle_A a$
\item[(ii)] $\langle \xi, \eta \rangle_A =  \langle \eta, \xi \rangle_A^*$
\item[(iii)]  $\langle \xi, \xi \rangle_A \geq 0$ and $\| \xi \| = \langle
\xi, \xi \rangle_A^{1/2}$
\end{enumerate}
for all $\xi, \eta \in X$ and $a \in A$.
For a Hilbert $A$-module $X$ we let $\Li(X)$ denote the
$C^*$-algebra of adjointable operators on $X$, and we let $\K (X)$ denote
the closed two-sided ideal of compact operators given by $$\K (X) := \clspan
\{ \Theta_{\xi,\eta}^X : \xi, \eta \in X \}$$ where $\Theta_{\xi,\eta}^X$ is
defined by $\Theta_{\xi,\eta}^X (\zeta) := \xi \langle \eta, \zeta
\rangle_A$.  When no confusion arises we shall often omit the superscript and
write $\Theta_{\xi,\eta}$ in place of $\Theta_{\xi,\eta}^X$.
\end{definition}

\begin{definition}
If $A$ is a $C^*$-algebra, then a \emph{$C^*$-correspondence} is a right
Hilbert $A$-module $X$ together with a $*$-homomorphism $\phi : A \to
\Li(X)$.  We consider $\phi$ as giving a left action of $A$ on $X$ by setting
$a \cdot x := \phi(a) x$.
\end{definition}

\begin{definition}
If $X$ is a $C^*$-correspondence over $A$, then a \emph{representation} of
$X$ into a $C^*$-algebra $B$ is a pair $(\pi, t)$ consisting of a
$*$-homomorphism $\pi : A \to B$ and a linear map $t : X \to B$ satisfying
\begin{enumerate}
\item[(i)] $t(\xi)^* t(\eta) = \pi(\langle \xi, \eta \rangle_A)$
\item[(ii)] $t(\phi(a)\xi) = \pi(a) t(\xi)$
\item[(iii)] $t(\xi a) = t(\xi) \pi(a)$
\end{enumerate}
for all $\xi, \eta \in X$ and $a \in A$.

Note that Condition (iii) follows from Condition (i) due to the equation $$\|
t(\xi) \pi(a) - t(\xi a) \|^2 = \| (t(\xi) \pi(a) - t(\xi a) )^* (t(\xi)
\pi(a) - t(\xi a)) \| = 0.$$  If $(\pi, t)$ is a representation of $X$ into
a $C^*$-algebra $B$, we let $C^*(\pi,t)$ denote the $C^*$-subalgebra of
$B$ generated by $\pi(A) \cup t(X)$.

A representation $(\pi, t)$ is said to be \emph{injective} if $\pi$ is
injective.  Note that in this case $t$ will also be isometric since $$\|
t(\xi) \|^2 = \| t(\xi)^* t(\xi) \| = \| \pi (\langle \xi, \xi \rangle_A) \| =
\| \langle \xi, \xi \rangle_A \| = \| \xi \|^2.$$

When $(\pi, t)$ is a representation of $X$ into $\B (\Hi)$ for a Hilbert
space $\Hi$, we say that $(\pi,t)$ is \emph{a representation of $X$ on $\Hi$}.
\end{definition}

In the literature a representation $(\pi,t)$ is sometimes referred to as a \emph{Toeplitz
representation} (See, e.g., \cite{FR} and \cite{FMR}.), and as an \emph{isometric
representation} \cite{MS}. However, here, all representations considered will be at least Toeplitz or isometric and so we drop the additional adjective.  We note that in \cite{FR} the authors show that given a correspondence $X$ over a $C^*$-algebra $A$, there is a $C^*$-algebra, denoted $\TX$ and a representation $({\pi}_X,t_X)$ of $X$ in $\TX$ that is universal in the following sense: $\TX$ is generated as a $C^*$-algebra by the ranges of ${\pi}_X$ and $t_X$, and given any representation $(\pi,t)$ in a $C^*$-algebra $B$, then there is a $C^*$-homomorphism of $\TX$ into $B$, denoted ${\rho}_{(\pi,t)}$, that is unique up to an inner automorphism of $B$, such that $\pi={{\rho}_{(\pi,t)}}\circ{\pi}_X$ and $t={{\rho}_{(\pi,t)}}\circ{t_X}$.  The $C^*$-algebra $\TX$ and the representation $({\pi}_X,t_X)$ are unique up to an obvious notion of isomorphism.  We call $\TX$ \emph{the} Toeplitz algebra of the correspondence $X$, but we call $({\pi}_X,t_X)$ \emph{a universal representation} of $X$ in $\TX$, with emphasis on the indefinite article, because at times we want to consider more than one.

\begin{definition}
For a representation $(\pi, t)$ of a $C^*$-correspondence $X$ on $B$ there
exists a $*$-homomorphism $\pi^{(1)} : \K (X) \to B$ with the property that
$$\pi^{(1)} (\Theta_{\xi,\eta}) = t(\xi) t(\eta)^*.$$  See \cite[p.~202]{Pim},
\cite[Lemma~2.2]{KPW}, and \cite[Remark~1.7]{FR} for details on the existence
of this $*$-homomorphism.  Also note that if $(\pi,t)$ is an injective
representation, then $\pi^{(1)}$ will be injective as well.
\end{definition}

\begin{definition}
For an ideal $I$ in a $C^*$-algebra $A$ we define $$I^\perp := \{ a \in A :
ab=0 \text{ for all } b \in I \}.$$  If $X$ is a $C^*$-correspondence over
$A$, we define an ideal $J(X)$ of $A$ by $J(X) := \phi^{-1}(\K(X))$.  We also
define an ideal $J_X$ of $A$ by $$J_X := J(X) \cap (\ker \phi)^\perp.$$  Note
that $J_X = J(X)$ when $\phi$ is injective, and that $J_X$ is the maximal
ideal on which the restriction of $\phi$ is an injection into $\K(X)$.
\end{definition}

\begin{definition}
If $X$ is a $C^*$-correspondence over $A$ and $K$ is an ideal in $J(X)$, then
we say that a representation $(\pi, t)$ is \emph{coisometric on $K$}, or is \emph{$K$-coisometric} if
$$\pi^{(1)} (\phi(a)) = \pi(a) \qquad \text{ for all $a \in K$}.$$
\end{definition}

In \cite[Proposition~1.3]{FMR} the authors show that given a correspondence $X$ over a $C^*$-algebra $A$, and an ideal $K$ of $A$ contained in $J(X)$, there is a $C^*$-algebra, denoted $\OKX$, and a representation $({\pi}_X,t_X)$ of $X$ in $\OKX$ that is coisometric on $K$ and is universal with this property, in the following sense: $\OKX$ is generated as a $C^*$-algebra by the ranges of ${\pi}_X$ and $t_X$, and given any representation $(\pi,t)$ of $X$ in a $C^*$-algebra $B$ that is $K$-coisometric, then there is a $C^*$-homomorphism of $\OKX$ into $B$, denoted ${\rho}_{(\pi,t)}$, that is unique up to an inner automorphism of $B$, such that $\pi={{\rho}_{(\pi,t)}}\circ{\pi}_X$ and $t={{\rho}_{(\pi,t)}}\circ{t_X}$.

\begin{definition}
The algebra $\OKX$, associated with an ideal $K$ in $J(X)$, is called the \emph{relative Cuntz-Pimsner algebra} determined by $X$ and the ideal $K$. Further, a representation $({\pi}_X,t_X)$ that is coisometric on $K$ and has the universal property just described is called a \emph{universal} $K$-coisometric representation of $X$.
\end{definition}

\begin{remark}
When the ideal $K$ is the zero ideal in $J(X)$, then the algebra $\OKX$ becomes $\TX$ and a universal $0$-coisometric representation of $X$ is simply a representation of $X$. Furthermore, if $X$ is a $C^*$-correspondence in which $\phi$ is injective,
then $\OX := \mathcal{O}(J_X,X)$ is precisely the augmented Cuntz-Pimsner
algebra of $X$ defined in \cite{Pim}. If $X$ is full, i.e., if $\clspan \{
\langle \xi,\eta \rangle_A : \xi, \eta \in X \} = A$, then the augmented
Cuntz-Pimsner algebra of $X$ and the Cuntz-Pimsner algebra of $X$ coincide.
Thus $\OX$ coincides with the Cuntz-Pimsner algebra of \cite{Pim} when $\phi$
is injective and $X$ is full. Whether or not $\phi$ is injective, a universal $J(X)$-coisometric representation is sometimes called a universal Cuntz-Pimsner covariant representation \cite[Definition~1.1]{FMR}.
\end{remark}

\begin{remark}
If $\mathcal{O}(K,X)$ is a relative Cuntz-Pimsner algebra associated to a
$C^*$-correspondence $X$, and if $(\pi, t)$ is a universal $K$-coisometric representation of
$X$, then for any $z \in \mathbb{T}$
$(\pi,zt)$ is also a universal $K$-coisometric representation. Hence by the universal property, there exists a homomorphism $\gamma_z
: \mathcal{O}(K,X) \rightarrow \mathcal{O}(K,X)$ such that
$\gamma_z(\pi(a))=\pi(a)$ for all $a \in A$ and $\gamma_z(t(\xi))=zt(\xi)$
for all $\xi \in X$.  Since $\gamma_{z^{-1}}$ is an inverse for this
homomorphism, we see that $\gamma_z$ is an automorphism. Thus we have an action
$\gamma : \mathbb{T} \rightarrow \operatorname{Aut} \mathcal{O}(K,X)$ with the
property that $\gamma_z (\pi(a)) = \pi(a)$ and $\gamma_z(t(\xi)) =
zt(\xi)$. Furthermore, a routine $\epsilon / 3$ argument shows that $\gamma$ is
strongly continuous. We call $\gamma$ the \emph{gauge action on
$\mathcal{O}(K,X)$}.
\end{remark}

%%%%%%%%%%%%%%%%%%%%%%%%%%%%%%%%%%%%%%%%%%%%%%%%%%%%%%%%%%
\section{Viewing graph $C^*$-algebras as Cuntz-Pimsner
algebras} \label{rel-gr-sec}
%%%%%%%%%%%%%%%%%%%%%%%%%%%%%%%%%%%%%%%%%%%%%%%%%%%%%%%%%%

Let $E := (E^0,E^1,r,s)$ be a directed graph with countable vertex set $E^0$,
countable edge set $E^1$, and range and source maps $r,s :E^1 \to E^0$. A
\emph{Cuntz-Krieger $E$-family} is a collection of partial isometries $\{
s_e : e \in E^1 \}$ with commuting range projections together with a
collection of mutually orthogonal projections $\{p_v : v \in E^0 \}$ that
satisfy

\begin{enumerate}
\item $s_e^*s_e=p_{r(e)}$ \text{ for all $e \in E^1$}
\item $s_es_e^* \leq p_{s(e)}$ \text{ for all $e \in E^1$}
\item $p_v = \sum_{\{e : s(e)=v \} } s_es_e^*$ \text{ for all $v \in E^0$
with $0 < | s^{-1}(v) | < \infty$}
\end{enumerate}

\noindent The \emph{graph algebra} $C^*(E)$ is the $C^*$-algebra generated by a
universal Cuntz-Krieger $E$-family (see \cite{KPRR,KPR,BPRS,FLR, BHRS}).

\begin{example}[The Graph $C^*$-correspondence]
If $E = (E^0,E^1,r,s)$ is a graph, we define $A:= C_0(E^0)$ and 
\begin{equation*}
X(E) := \{ x : E^1 \to \mathbb{C} : \text{ the function } v \mapsto \sum_{
\{f \in E^1: r(f) = v \} } |x(f)|^2 \text{ is in $C_0(E^0)$} \ \}.
\end{equation*}
Then $X(E)$ is a $C^*$-correspondence over $A$ with the operations 
\begin{align*}
(x \cdot a)(f) &:= x(f) a(r(f)) \text{ for $f \in E^1$} \\
\langle x, y \rangle_A(v) &:= \sum_{ \{ f \in E^1: r(f) = v \} }\overline{x(f)}y(f) \text{ for $f \in E^1$} \\
(a \cdot x)(f) &:= a(s(f)) x(f) \text{ for $f \in E^1$}
\end{align*}
and we call $X(E)$ the \emph{graph $C^*$-correspondence} associated to
$E$. Note that we could write $X(E) = \bigoplus_{v \in E^0}^0
\ell^2(r^{-1}(v))$ where this denotes the $C_0$ direct sum (sometimes called
the restricted sum) of the $\ell^2(r^{-1}(v))$'s. Also note that $X(E)$ and $A$
are spanned by the point masses $\{\delta_f : f \in E^1 \}$ and $\{ \delta_v :
v \in E^0 \}$, respectively.
\end{example}

\begin{theorem}[{\text{\protect\cite[Proposition~12]{FLR}}}]
If $E$ is a graph with no sinks, and $X(E)$ is the associated graph
$C^*$-correspondence, then $\mathcal{O}(J(X(E)), X(E)) \cong C^*(E)$.
Furthermore, if $(\pi_X, t_X)$ is a universal $J(X(E))$-coisometric representation, then $\{ t_X(\delta_e), \pi_X(\delta_v) \}$ is a universal
Cuntz-Krieger $E$-family in $\mathcal{O}(J(X(E)),X(E))$.
\end{theorem}

\noindent It was shown in \cite[Proposition~4.4]{FR} that 
\begin{equation*}
J(X(E)) = \operatorname{\overline{\mathrm{span}}} \{ \delta_v : |s^{-1}(v)| <
\infty \}
\end{equation*}
and if $v$ emits finitely many edges, then 
\begin{equation*}
\phi(\delta_v) = \sum_{\{f \in E^1 : s(f) =v\} } \Theta_{\delta_f,\delta_f}
\quad \text{ and } \quad \pi_X(\phi(\delta_v)) = \sum_{\{f \in E^1 : s(f) =v\}
} t_X(\delta_f)t_X(\delta_f)^*.
\end{equation*}
Furthermore, one can see that $\delta_v \in \ker \phi$ if and only if $v$ is
a sink in $E$. Also $\delta_v \in \operatorname{\overline{\mathrm{span}}} \{ \langle
x, y \rangle_A \}$ if and only if $v$ is a source, and since $\delta_{s(f)}
\cdot \delta_f = \delta_f$ we see that $\operatorname{\overline{\mathrm{span}}} A
\cdot X =X$ and $X(E)$ is essential. These observations show that we have
the following correspondences between the properties of the graph $E$ and
the properties of the graph $C^*$-correspondence $X(E)$.

\begin{center}
\begin{tabular}{|c|c|}
\hline
\textbf{Property of $\boldsymbol{X(E)}$} &
\textbf{Property of $\boldsymbol{E}$} \\ \hline
$\phi(\delta_v) \in \mathcal{K}(X(E))$ & $v$ emits a finite number of edges
\\ 
$\phi(A) \subseteq \mathcal{K}(X(E))$ & $E$ is row-finite \\ 
$\phi$ is injective & $E$ has no sinks \\ 
$X(E)$ is full  & $E$ has no sources \\ 
$X(E)$ is essential & always \\ \hline
\end{tabular}
\end{center}

\begin{remark}
If $E$ is a graph with no sinks, then $\mathcal{O}(J(X(E)), X(E))$ is
canonically isomorphic to $C^*(E)$. When $E$ has sinks, this will not be the
case.  If $(\pi, t)$ is the universal $J(X(E))$-coisometric representation of $X(E)$, then it will be the case that $\{ t(\delta_e),
\pi (\delta_v) \}$ is a Cuntz-Krieger $E$-family.  However, when $v$ is a
sink in $E$, $\phi(\delta_v) = 0$ and thus $\pi (\delta_v) =
\pi^{(1)} (\phi(\delta_v)) = 0$.  Consequently, $\{ t(\delta_e),
\pi (\delta_v) \}$ will not be a universal Cuntz-Krieger $E$-family when $E$
has sinks.

However, if $E$ is a graph with sinks, then we see that $\phi(\delta_v) = 0$ if
and only if $v$ is a sink, and $\delta_v \in (\ker \phi)^\perp$ if and only if
$v$ is not a sink.  Thus $$J_{X(E)} = \operatorname{\overline{\mathrm{span}}}\{
\delta_v : 0 < | s^{-1}(v) | < \infty
\}$$ and a proof similar to that in \cite[Proposition~4.4]{FR} shows that
$\mathcal{O}_{X(E)} := \mathcal{O} (J_{X(E)}, X(E)) $ is isomorphic to
$C^*(E)$.  Furthermore, if $(\pi_X, t_X)$ is a universal $J(X(E))$-coisometric representation of
$X(E)$, then $\{t_X(\delta_e),
\pi_X(\delta_v) \}$ is a universal Cuntz-Krieger $E$-family in
$\mathcal{O}_{X(E)}$.
\end{remark}

%-------------------------------------------
\subsection{Relative Graph Algebras} \label{rel-gr-subsec}
%-------------------------------------------

We shall now examine relative Cuntz-Pimsner algebras in the context of
graph algebras. If $E$ is a graph and $X(E)$ is the associated graph
$C^*$-correspondence, then $J_{X(E)} := \operatorname{\overline{\mathrm{span}}}
\{ \delta_v :  0 < |s^{-1}(v)| < \infty \}$. If $K$ is an ideal in $J_{X(E)}$,
then $K = \operatorname{\overline{\mathrm{span}}} \{ \delta_v : v \in V \}$ for
some subset $V$ of vertices which emit a finite and nonzero number of edges. If
$(\mathcal{O}(K,X(E)), t_X, \pi_X)$ is the relative Cuntz-Pimsner algebra
determined by $K$, then the relation
$\pi_X(\delta_v) = \sum_{s(e) =v} t_X(\delta_e)t_X(\delta_e)^*$ will hold only
for vertices $v \in V$. This motivates the following definition.

\begin{definition}
Let $E = (E^0,E^1,r,s)$ be a graph and define $R(E) := \{ v \in E^0 : 0 < |
s^{-1}(v) | < \infty \}$. For any $V \subseteq R(E)$ we define a \emph{Cuntz-Krieger
$(E,V)$-family} to be a collection of mutually orthogonal projections $\{p_v : v \in
E^0 \}$ together with a collection of partial isometries $\{ s_e : e \in E^1 \}$ that
satisfy

\begin{enumerate}
\item $s_e^*s_e = p_{r(e)}$ \text{ for $e \in E^1$}

\item $s_es_e^* < p_{s(e)}$ \text{ for $e \in E^1$}

\item $p_v = \sum_{s(e)=v} s_es_e^*$ \text{ for all $v \in V$}
\end{enumerate}

We refer to a Cuntz-Krieger $(E,R(E))$-family as simply a Cuntz-Krieger
$E$-family, and we refer to a Cuntz-Krieger $(E, \emptyset)$-family as a
Toeplitz-Cuntz-Krieger family.
\end{definition}

\begin{definition}
\label{rel-ga} If $E$ is a graph and $V \subseteq R(E)$, then we define the 
\emph{relative graph algebra} $C^*(E,V)$ to be the $C^*$-algebra generated
by a universal Cuntz-Krieger $(E,V)$-family.
\end{definition}

The existence of $C^*(E,V)$ can be proven by adapting the argument for the
existence of graph algebras in \cite{KPR}, or by realizing $C^*(E,V)$ as a
relative Cuntz-Pimsner algebra.

Note that $C^*(E,R(E))$ is the graph algebra $C^*(E)$, and $C^*(E,\emptyset)$
is the Toeplitz algebra defined in \cite[Theorem~4.1]{FR} (but different
from the Toeplitz algebra defined in \cite{EL}). It is also the case that if
$\{s_e, p_v\}$ is a universal Cuntz-Krieger $(E,V)$-family, then whenever $v
\in R(E) \backslash V$ we have $p_v > \sum_{s(e)=v} s_es_e^*$.

\begin{definition}
Let $E = (E^0, E^1,r,s)$ be a graph and $V \subseteq R(E)$. We define the
graph $E_V$ to be the graph with vertex set $E_V^0 := E^0 \cup \{
v^{\prime}: v \in R(E) \backslash V \}$, edge set $E^1 \cup \{e^{\prime}: e
\in E^1 \text{ and } r(e) \in R(E) \backslash V \}$, and $r$ and $s$
extended to $E_V^1$ by defining $s(e^{\prime}) := s(e)$ and $r(e^{\prime})
:= r(e)^{\prime}$.
\end{definition}

Roughly speaking, when forming $E_V$ one takes $E$ and adds a sink for each
element $v \in R(E) \backslash V$ as well as edges to this sink from each
vertex that feeds into $v$.

\begin{theorem}
\label{rel-gas-are-gas} If $E$ is a graph and $V \subseteq R(E)$, then the
relative graph algebra $C^*(E,V)$ is canonically isomorphic to the graph
algebra $C^*(E_V)$.
\end{theorem}

\begin{proof}
Let $\{s_e, p_v: e \in E^1, v \in E^0 \}$ be a generating Cuntz-Krieger
$(E,V)$-family in $C^*(E,V)$. For $w \in E_V^0$ and $f \in E_V^1$ define 
\begin{align*}
q_w &:= 
\begin{cases}
p_v & \text{ if $w \notin R(E) \backslash V$} \\ 
\sum_{ \{e \in E^1 : s(e) = w \} } s_es_e^* & \text{ if $w \in R(E)
\backslash V$} \\ 
p_v - \sum_{ \{ e \in E^1 : s(e) = v \} } s_es_e^* & \text{ if $w=v^{\prime}$
for some $v \in R(E) \backslash V$.}
\end{cases}
\\
t_f &:= 
\begin{cases}
s_f q_{r(f)} & \text{ if $f \in E^1$} \\ 
s_eq_{r(e)^{\prime}} & \text{ if $f = e^{\prime}$ for some $e \in E^1$.}
\end{cases}
\end{align*}
It is straightforward to check that $\{ t_f, q_w : f \in E_V^1, w \in E_V^0
\}$ is a Cuntz-Krieger $E_V$-family in $C^*(E,V)$. Thus by the universal
property there exists a homomorphism $\alpha : C^*(E_V) \to C^*(E,V)$ taking
the generators of $C^*(E_V)$ to $\{ t_f, q_w\}$. By the gauge-invariant
uniqueness theorem \cite[Theorem~2.1]{BHRS} $\alpha$ is injective.
Furthermore, whenever $v \in R(E) \backslash V$ we see that $p_v = q_v +
q_{v^{\prime}}$ and whenever $r(e) \in R(E) \backslash V$ we see that $s_e =
t_e + t_{e^{\prime}}$. Thus $\{ q_w, t_f \}$ generates $C^*(E,V)$ and $\alpha
$ is surjective. Consequently $\alpha$ is an isomorphism.
\end{proof}

This theorem shows that the class of relative graph algebras is the same as
the class of graph algebras. Thus we gain no new $C^*$-algebras by
considering relative graph algebras in place of graph algebras. However, we
maintain that relative graph algebras are still useful and arise naturally
in the study of graph algebras. In particular, we give three examples of
common situations in which relative graph algebras prove convenient.

\begin{example}[Subalgebras of Graph Algebras]
Let $E = (E^0,E^1,r,s)$ be a graph and let $\{s_e,p_v : e\in E^1, v \in E^0
\}$ be a generating Cuntz-Krieger $E$-family in $C^*(E)$. If $F=(F^0,F^1,r_F,s_F)$ is
a subgraph of $E$, and $A$ denotes the $C^*$-subalgebra of $C^*(E)$ generated by $\{
s_e, p_v : e \in F^1, v \in F^0 \}$, then it is well-known that $A$ is a graph
algebra (but not necessarily the $C^*$-algebra associated to $F$). In fact, we
see that for any $v \in F^0$, the sum $\sum_{\{ e \in F^1 : s_F(e) = v \}
}s_es_e^*$ may not add up to $p_v
$ because some of the edges in $s^{-1}(v)$ may not be in $F$. However, if we
let $V := \{ v \in R(F) : s_F^{-1}(v) = s^{-1}(v) \}$. Then $\{s_e,p_v :
e\in F^1, v \in F^0 \}$ is a Cuntz-Krieger $(F,V)$-family and $A \cong
C^*(F,V)$.

These subalgebras arise often in the study of graph algebras. In
\cite[Lemma~2.4]{HS3} they were realized as graph algebras by the method shown in the
proof of Theorem~\ref{rel-gas-are-gas}, and in \cite[Lemma~1.2]{RS} these subalgebras
were realized as graph algebras by using the notion of a dual graph. In both of these
instances it would have been convenient to have used relative graph algebras.
Realizing the subalgebra as $C^*(F,V)$ would have provided an economy of notation as
well as a more direct analysis of the subalgebras under consideration.
\end{example}

\begin{example}[Spielberg's Toeplitz Graph Algebras]
In \cite{Spi2} Spielberg introduced a notion of a Toeplitz graph groupoid
and a Toeplitz graph algebra. The Toeplitz graph algebras defined in
\cite[Definition~2.17]{Spi2} are relative graph algebras as defined in
Definition~\ref{rel-ga} (see \cite[Theorem~2.9]{Spi2}). Spielberg also made use of his
Toeplitz graph algebras in \cite{Spi} to construct graph algebras with a specified
$K$-theory.
\end{example}

\begin{example}[Quotients of Graph Algebras]
If $E=(E^0,E^1,r,s)$ is a row-finite graph and $H$ is a saturated hereditary
subset of vertices of $E$, then it follows from \cite[Theorem~4.1(b)]{BPRS}
that $C^*(E) / I_H \cong C^*(F)$ where $F$ is the subgraph defined by 
\begin{equation*}
F^0 := E^0 \backslash H \qquad \qquad F^1 := \{ e \in E^1 : r(e) \notin H \}.
\end{equation*}
If $E$ is not row-finite, then this is not necessarily the case. The
obstruction is due to the vertices in the set 
\begin{equation*}
B_H := \{v \in E^0 \,|\, 
\hbox{$v$ is an infinite-emitter and $0 < |s^{-1}(v)
\cap r^{-1}(E^0 \setminus H)| < \infty$}\}.
\end{equation*}
In fact, if $\{s_e,p_v \}$ is a generating Cuntz-Krieger $E$-family in $C^*(E)$, then
the cosets $\{s_e + I_H, p_v + I_H : v \notin H, r(e) \notin H
\}$ will have the property that $p_v + I_H \geq \sum_{e \in E \backslash H :
s(e)= v \} } (s_e+I_H) (s_e+I_H)^*$ with equality occurring if and only if $v \in R(F)
\backslash B_H$. Thus it turns out that $\{s_e + I_H, p_v + I_H : v \notin H, r(e)
\notin H \}$ will be a Cuntz-Krieger $(F, R(F) \backslash B_H)$-family and $C^*(E) /
I_H \cong C^*(F, R(F) \backslash B_H)$.

The quotient $C^*(E)/I_H$ was realized as a graph algebra
in \cite[Proposition~3.4]{BHRS} by a technique similar to that used in the proof of
Theorem~\ref{rel-gas-are-gas}. However, relative graph algebras provide a more natural
context for describing these quotients.
\end{example}

In addition to their applications in the situations mentioned above,
relative graph algebras can be useful for another reason. Since any relative
graph algebra is canonically isomorphic to a graph algebra, we see that for
every theorem about graph algebras there will be a corresponding theorem for
relative graph algebras. Thus the relative graph algebras provide a class of
relative Cuntz-Pimsner algebras that are well understood.  With this in mind, we
shall now state a version of the Gauge-Invariant Uniqueness Theorem for
relative graph algebras.

\begin{theorem}[Gauge-Invariant Uniqueness for Relative Graph Algebras]
\label{GIU-rel-ga} Let $E=(E^0,E^1,r,s)$ be a graph and $V \subseteq R(E)$.
Also let $\{s_e,p_v : e \in E^1, v \in E^0 \}$ and let $\gamma : \mathbb{T}
\to \operatorname{Aut} C^*(E,V)$ denote the gauge action on $C^*(E,V)$. If $\rho :
C^*(E,V) \to A$ is a $*$-homomorphism between $C^*$-algebras that satisfies
\begin{enumerate}
\item $\rho(p_v) \neq 0$ for all $v \in E^0$
\item $\rho( p_v - \sum_{s(e)=v} s_es_e^* ) \neq 0$ for all $v \in V$
\item there exists a strongly continuous action $\beta : \mathbb{T} \to
\operatorname{Aut} A$ such that $\beta_z \circ \rho = \rho \circ \gamma_z$ for
all $z \in \mathbb{T}$.
\end{enumerate}
then $\rho$ is injective.
\end{theorem}

\begin{proof}
By Theorem~\ref{rel-gas-are-gas} there exists an isomorphism $\alpha :C^*(E_V)
\to C^*(E,V)$ and a generating Cuntz-Krieger $E_V$-family $\{t_e,q_w\}$ for
which 
\begin{align*}
\alpha(q_w) &:= 
\begin{cases}
p_v & \text{ if $w \notin R(E) \backslash V$} \\ 
\sum_{ \{e \in E^1 : s(e) = w \} } s_es_e^* & \text{ if $w \in R(E)
\backslash V$} \\ 
p_v - \sum_{ \{ e \in E^1 : s(e) = v \} } s_es_e^* & \text{ if $w=v^{\prime}$
for some $v \in R(E) \backslash V$.}
\end{cases}
\\
\alpha(t_f) &:= 
\begin{cases}
s_f q_{r(f)} & \text{ if $f \in E^1$} \\ 
s_eq_{r(e)^{\prime}} & \text{ if $f = e^{\prime}$ for some $e \in E^1$.}
\end{cases}
\end{align*}
To show that $\rho$ is injective, it suffices to show that $\rho \circ \alpha$
is injective. We shall do this by applying the gauge-invariant uniqueness
theorem for graph algebras \cite[Theorem~2.1]{BHRS} to $\rho \circ \alpha$. Now
clearly if $w \notin R(E) \backslash V$, then $\rho \circ \alpha (q_w)
\neq 0$ by $(1)$. If $w = v^{\prime}$, then $\rho \circ \alpha (q_w) \neq 0$
by $(2)$. Furthermore, if $w \in R(E) \backslash V$ then $\rho \circ \alpha
(q_w) = 0$ implies that $\rho(\sum_{s(e)=w} s_es_e^*) = 0$ and thus for any $f \in
s^{-1}(v)$ we have 
\begin{equation*}
\rho(s_f) = \rho(\sum_{s(e)=w} s_es_e^*)\rho(s_f) = 0.
\end{equation*}
But then $\rho(p_{r(f)}) = \rho(s_f^*s_f)=0$ which contradicts $(1)$. Hence
we must have $\rho \circ \alpha (q_w) \neq 0$. Finally, if $\gamma^{\prime}$
denotes the gauge action on $C^*(E_V)$, then by checking on generators we
see that $\beta_z \circ (\rho \circ \alpha) = (\rho \circ \alpha) \circ
\gamma^{\prime}_z$. Therefore, $\rho \circ \alpha$ is injective by the gauge
invariant uniqueness theorem for graph algebras, and consequently $\rho$ is
injective.
\end{proof}

We have shown in Theorem~\ref{rel-gas-are-gas} that every relative graph
algebra is isomorphic to a graph algebra.  More generally, Katsura has shown
in \cite{Kat5} that every relative Cuntz-Pimsner algebra is isomorphic to the
$C^*$-algebra associated to a correspondence; that is, if $\mathcal{O}(K,X)$ is
a relative Cuntz-Pimsner algebra, then there exists a $C^*$-correspondence $X'$
such that $\mathcal{O}_{X'} := \mathcal{O}(J_{X'},X')$ is isomorphic to
$\mathcal{O}(K,X)$.  In Theorem~\ref{GIU} we shall prove a gauge-invariant
uniqueness theorem for $C^*$-algebras associated to correspondences. 
Afterwards, in Remark~\ref{rel-GIU}, we shall use Katsura's analysis in
\cite{Kat5} to give an interpretation of Theorem~\ref{GIU-rel-ga} in
the context of relative Cuntz-Pimsner algebras.

%%%%%%%%%%%%%%%%%%%%%%%%%%%%%%%%%%%%%%%%%%%%%%%%%%%%%%%%%%%%%%%
\section{Adding Tails to $C^*$-correspondences} \label{add-tails-sec}
%%%%%%%%%%%%%%%%%%%%%%%%%%%%%%%%%%%%%%%%%%%%%%%%%%%%%%%%%%%%%%%

If $E$ is a graph and $v$ is a vertex of $E$, then by \emph{adding a tail to
$v$} we mean attaching a graph of the form
\begin{equation*}
\xymatrix{ v \ar[r]^{e_1} & v_1 \ar[r]^{e_2} & v_2 \ar[r]^{e_3} & v_3
\ar[r]^{e_4} & \cdots\\ }
\end{equation*}
to $E$.  It was shown in \cite[\S1]{BPRS} that if $F$ is the graph formed by
adding a tail to every sink of $E$, then $F$ is a graph with no sinks and
$C^*(E)$ is canonically isomorphic to a full corner of $C^*(F)$.  The
technique of adding tails to sinks is a simple but powerful tool in the
analysis of graph algebras.  In the proofs of many results it allows one to
reduce to the case in which the graph has no sinks and thereby avoid certain
complications and technicalities.

Our goal in this section is to develop a process of ``adding tails to sinks"
for $C^*$-correspondences, so that given any $C^*$-correspondence $X$ we may
form a $C^*$-correspondence $Y$ with the property that the left action of
$Y$ is injective and $\OX$ is canonically isomorphic to a full corner in
$\OY$.

\begin{definition}
Let $X$ be a $C^*$-correspondence over $A$ with left action $\phi : A \to
\Li(X)$, and let $I$ be an ideal in $A$.  We define \emph{the tail determined
by $I$} to be the $C^*$-algebra 
\begin{equation*}
T := I^{(\mathbb{N})}
\end{equation*}
where $I^{(\mathbb{N})}$ denotes the $c_0$-direct sum of countably
many copies of the ideal $I$.  We shall denote the elements of $T$ by 
\begin{equation*}
\vec{f} := (f_1, f_2, f_3, \ldots )
\end{equation*}
where each $f_i$ is an element of $I$. We shall consider $T$ as a
right Hilbert $C^*$-module over itself (see \cite[Example~2.10]{RW}).  We define
$Y := X \oplus T$ and $B := A \oplus T$. Then $Y$ is a right Hilbert $B$-module
in the usual way; that is, the right action is given by 
\begin{equation*}
(\xi,\vec{f}) \cdot (a, \vec{g}) := (\xi \cdot a, \vec{f}\vec{g}) \qquad 
\text{for $\xi \in X$, $a \in A$, and $\vec{f},\vec{g} \in T$}
\end{equation*}
and the inner product is given by 
\begin{equation*}
\langle (\xi, \vec{f}), (\nu, \vec{g}) \rangle_B := (\langle \xi, \nu
\rangle_A, \vec{f}^*\vec{g} ) \qquad \text{for $\xi, \nu \in X$ and
$\vec{f},\vec{g} \in T$.}
\end{equation*}
Furthermore, we shall make $Y$ into a $C^*$-correspondence over $B$ by
defining a left action $\phi_B : B \to \mathcal{L}(Y)$ as 
\begin{equation*}
\phi_B(a, \vec{f}) (\xi, \vec{g}) := (\phi(a)(\xi), (a g_1, 
f_1 g_2, f_2 g_3, \ldots )) \ \text{ for $a
\in A$, $\xi \in X$, and $\vec{f}, \vec{g} \in T$.}
\end{equation*}
We call $Y$ \emph{the $C^*$-correspondence formed by adding the tail $T$ to
$X$}.
\end{definition}

\begin{lemma} \label{phi-B-inj}
Let $X$ be a $C^*$-correspondence over $A$, and let $T :=(\ker \phi)^{(\N)}$
be the tail determined by $\ker \phi$.  If $Y:= X \oplus T$ is the
$C^*$-correspondence over $B := A \oplus T$ formed by adding the tail $T$ to
$X$, then the left action $\phi_B : B \to \Li(Y)$ is injective. 
Consequently, $J_Y = J(Y)$ and $\mathcal{O}_Y = \mathcal{O}(J(Y),Y)$ is equal
to the $C^*$-algebra defined by Pimsner in \cite{Pim}.
\end{lemma}

\begin{proof}
If $(a, \vec{f}) \in \ker \phi_B$, then for all $\xi \in X$ we have
$$(\phi(a) \xi, \vec{0}) = \phi_B(a,\vec{f}) (\xi, \vec{0}) = (0, \vec{0})$$
so that $\phi(a) \xi = 0$ and $a \in \ker \phi$.  Thus $(0, (a, f_1, f_2,
\ldots)) \in X \oplus T$ and $$(0, (aa^*, f_1f_1^*, f_2f_2^*, \ldots)) =
\phi_B(a, \vec{f}) (0, (a^*, f_1, f_2, \ldots)) = (0, \vec{0})$$ so that $\|
a \|^2 = \| aa^*\| = 0$ and $\| f_i \|^2 = \| f_if_i^* \| = 0$ for all $i \in
\N$.  Consequently, $a = 0$ and $\vec{f} = \vec{0}$ so that $\phi_B$ is
injective.
\end{proof}

\begin{theorem} \label{pi-t-extend}
Let $X$ be a $C^*$-correspondence over $A$, and let $T :=(\ker \phi)^{(\N)}$
be the tail determined by $\ker \phi$.  Also let $Y:= X \oplus T$ be the
$C^*$-correspondence over $B := A \oplus T$ formed by adding the tail $T$ to
$X$.

\begin{enumerate}
\item[(a)] If $(\pi,t)$ is a $J_X$-coisometric representation of $X$ on a Hilbert
space $\Hi_X$, then there is a Hilbert space $\Hi_Y = \Hi_X \oplus \Hi_T$ and a $J(Y)$-coisometric representation $(\tilde{\pi}, \tilde{t})$ of $Y$ on $\Hi_Y$ with the property that $\tilde{\pi}|_X = \pi$ and $\tilde{t}|_A = t$.

\item[(b)]  If $(\tilde{\pi}, \tilde{t})$ is a $J(Y)$-coisometric representation of
$Y$ into a $C^*$-algebra $C$, then $(\tilde{\pi}|_A, \tilde{t}|_X)$ is a $J_X$-coisometric
representation of $X$ into $C$.  Furthermore, if $\tilde{\pi}|_A$ is injective, then $\tilde{\pi}$ is injective.

\item[(c)]  Let $(\pi_Y,t_Y)$ be a universal $J(Y)$-coisometric representation of $Y$.  Then $(\pi, t) := (\pi_Y|_A, t_Y|_X)$ is a $J_X$-coisometric representation of $X$ in $C^*(\pi_Y,t_Y)$.  Furthermore, $\rho_{(\pi,t)} : \mathcal{O}_X \to C^*(\pi_X,t_X) \subseteq \mathcal{O}_Y$ is an isomorphism onto the $C^*$-subalgebra of
$\mathcal{O}_Y$ generated by $$\{ \pi_Y( a, \vec{0}), t_Y(\xi, \vec{0}) : a \in
A \text{ and }  \xi \in X \}$$ and this $C^*$-subalgebra is a full corner of
$\mathcal{O}_Y$.  Consequently, $\mathcal{O}_X$ is naturally isomorphic to a
full corner of $\mathcal{O}_Y$.
\end{enumerate}

\end{theorem}

\begin{corollary} \label{univ-repn-inj}
If $X$ is a $C^*$-correspondence and $(\pi_X, t_X)$ is a universal $J(X)$-coisometric
representation of $X$, then $(\pi_X, t_X)$ is
injective.
\end{corollary}

\begin{proof}
By the theorem $(\pi_X, t_X)$ extends to a universal $J(Y)$-coisometric representation $(\pi_Y,
t_Y)$ of $Y$. Since $\phi_B$ is injective by
Lemma~\ref{phi-B-inj} it follows from \cite[Corollary~6.2]{FMR} that
$(\pi_Y, t_Y)$ is injective.  Consequently, $\pi_X = \pi_Y |_A$ is injective.
\end{proof}

\noindent To prove this theorem we shall need a number of lemmas.

\begin{lemma} \label{identify-J(Y)}
Let $X$ be a $C^*$-correspondence and let $T := (\ker \phi)^{(\N)}$ be the
tail determined by $\ker \phi$.  Also let $Y := X \oplus T$ be the
$C^*$-correspondence over $B:= A \oplus T$ formed by adding the tail $T$ to $X$.
Then for any $(a, \vec{f}) \in Y$ we have that  $(a, \vec{f}) \in J(Y)$ if and
only if $a = a_1 + a_2$ with $a_1 \in J_X$ and $a_2 \in \ker \phi$.
\end{lemma}

\begin{proof}
Suppose $a = a_1 + a_2$ with $a_1 \in J_X$ and $a_2 \in \ker \phi$.  Then we
may write $\phi(a_1) = \lim_n \sum_{k=1}^{N_n} \Theta^X_{\xi_{n,k}, \eta_{n,k}}$
for some $\xi_{n,k}, \eta_{n,k} \in X$.  But then $$\phi_B(a_1, \vec{0}) =
\lim_n \sum_{k=1}^{N_n} \Theta^Y_{(\xi_{n,k}, \vec{0}), (\eta_{n,k}, \vec{0})}
\in \K(Y).$$ In addition, since $a_2 \in \ker \phi$ we see that if we let $\{
\vec{e}_\lambda \}_{\lambda \in \Lambda}$ be an approximate unit for $T$ with
$\vec{e}_\lambda = (e_\lambda^1, e_\lambda^2, \ldots )$ for each $\lambda$,
then $$\phi_B(a_2 ,\vec{f}) =
\lim_\lambda \Theta^Y_{(0, (a_2, f_1, f_2, \ldots)), (0, (e_\lambda^1,
e_\lambda^2, \ldots))} \in \K(Y).$$  Thus $\phi_B(a,\vec{f})  = \phi_B(a_1,
\vec{0}) + \phi_B(a_2, \vec{f})) \in \K(Y)$.

Conversely, suppose that $\phi_B(a, \vec{f}) \in \K(Y)$.  Then we may write
$$\phi_B(a, \vec{f}) = \lim_n \sum_{k=1}^{N_n} \Theta_{(\xi_{n,k},
\vec{f}_{n,k}), (\eta_{n,k}, \vec{g}_{n,k})}^Y.$$  If we write $\vec{f}_{n,k} =
(f_{n,k}^1, f_{n,k}^2, \ldots )$ and $\vec{g}_{n,k} = (g_{n,k}^1, g_{n,k}^2,
\ldots )$ then for any $(\xi, \vec{g}) \in X \oplus T$ we have that
\begin{align}
(\phi(a) \xi, (ag_1, f_1g_2, \ldots)) 
& = \phi_B(a, \vec{f}) (\xi, \vec{g}) \notag \\ 
& = \lim_n \sum_{k=1}^{N_n} \Theta_{(\xi_{n,k}, \vec{f}_{n,k}), (\eta_{n,k},
\vec{g}_{n,k})}^Y (\xi, \vec{g}) \notag \\ 
& = \lim_n \sum_{k=1}^{N_n} ( (\xi_{n,k} \langle \eta_{n,k}, \xi
\rangle_A, (f_{n,k}^1{g_{n,k}^1}^* g_1, f_{n,k}^2{g_{n,k}^2}^* g_2, \ldots))
\notag \\
& = \lim_n \sum_{k=1}^{N_n} ( (\Theta_{\xi_{n,k},
\eta_{n,k}}^X \xi, (f_{n,k}^1{g_{n,k}^1}^* g_1, f_{n,k}^2{g_{n,k}^2}^* g_2,
\ldots)). \label{a2-mult}
\end{align}
Now since the operator norm on $\Li (Y)$ dominates the operator norm on
$\Li(X)$, we see that $\lim_n \sum_{k=1}^{N_n} \Theta_{\xi_{n,k},
\eta_{n,k}}^X$ converges and $\phi(a) = \lim_n \sum_{k=1}^{N_n}
\Theta_{\xi_{n,k}, \eta_{n,k}}^X$.  Thus $a \in \K(X)$.

Furthermore, if $\{e_\lambda\}_{\lambda \in \Lambda}$ is an approximate unit
for $\ker \phi$, then for any $n,m \in \N$ we have
\begin{align*}
& \Big\| ( \sum_{k=1}^{N_n} f_{n,k}^1{g_{n,k}^1}^* -
\sum_{k=1}^{N_m} f_{m,k}^1{g_{m,k}^1}^* ) e_\lambda \Big\| \\ 
= \ & \Big\| \big(\sum_{k=1}^{N_n} \Theta_{(\xi_{n,k}, \vec{f}_{n,k}),
(\eta_{n,k}, \vec{g}_{n,k})}^Y - \sum_{k=1}^{N_m}
\Theta_{(\xi_{m,k}, \vec{f}_{m,k}), (\eta_{m,k},
\vec{g}_{m,k})}^Y \big) (0, (e_\lambda, 0, 0, \ldots)) \Big\| \\
\leq \ & \Big\| \sum_{k=1}^{N_n} \Theta_{(\xi_{n,k}, \vec{f}_{n,k}),
(\eta_{n,k}, \vec{g}_{n,k})}^Y - \sum_{k=1}^{N_m}
\Theta_{(\xi_{m,k}, \vec{f}_{m,k}), (\eta_{m,k},
\vec{g}_{m,k})}^Y \Big\| \ \Big\| (0, (e_\lambda, 0, 0, \ldots)) \Big\| \\
= \ & \Big\| \sum_{k=1}^{N_n} \Theta_{(\xi_{n,k}, \vec{f}_{n,k}), (\eta_{n,k},
\vec{g}_{n,k})}^Y - \sum_{k=1}^{N_m}
\Theta_{(\xi_{m,k}, \vec{f}_{m,k}), (\eta_{m,k}, \vec{g}_{m,k})}^Y \Big\|
\end{align*} 
for all $\lambda \in \Lambda$.  Taking the limit with respect to $\lambda$
shows that $$\big\| \sum_{k=1}^{N_n} f_{n,k}^1{g_{n,k}^1}^* - \sum_{k=1}^{N_m}
f_{m,k}^1{g_{m,k}^1}^* \big\| \leq 
\big\| \sum_{k=1}^{N_n} \Theta_{(\xi_{n,k}, \vec{f}_{n,k}), (\eta_{n,k},
\vec{g}_{n,k})}^Y - \sum_{k=1}^{N_m}
\Theta_{(\xi_{m,k}, \vec{f}_{m,k}), (\eta_{m,k}, \vec{g}_{m,k})}^Y) \big\|.$$ 
Since the $\sum_{k=1}^{N_n} \Theta_{(\xi_{n,k}, \vec{f}_{n,k}), (\eta_{n,k},
\vec{g}_{n,k})}^Y$'s converge in the operator norm on $\Li(Y)$, this inequality
implies that $\sum_{k=1}^{N_n} f_{n,k}^1{g_{n,k}^1}^*$ converges to an element
in $\ker \phi$.  If we let $a_2 = \lim_n \sum_{k=1}^{N_n} f_{n,k}^1{g_{n,k}^1}^*
\in \ker \phi$, then Eq.(\ref{a2-mult}) shows that $ag = a_2g$ for all $g \in
\ker \phi$.  But then $a_1 := a-a_2 \in (\ker \phi)^\perp$, and consequently
$a_1 \in J_X$.  Since $a = a_1 + a_2$ the proof is complete. 
\end{proof}

\begin{lemma}  \label{b-inj-lemma}
Let $(\tilde{\pi}, \tilde{t})$ be a representation of $Y$ which is
coisometric on $\ker \phi \oplus T$, and suppose that
$\tilde{\pi}|_A$ is injective.  For any $f \in \ker \phi$ we
define $\epsilon_i(f) := (0, \ldots, 0, f, 0, \ldots) \in T$ where $f$ appears
in the $i$\textsuperscript{th} position.  Then for every $i \in \N$ and for
every $f \in  \ker \phi$, the equation $\tilde{\pi}(0, \epsilon_i(f)) = 0$
implies that $f = 0$.
\end{lemma}

\begin{proof}
First note that it suffices to prove the lemma for $f \geq 0$, because if
$\tilde{\pi}(0, \epsilon_i(f)) = 0$ then $\tilde{\pi}(0, \epsilon_i(f f^*)) =
\tilde{\pi}(0, \epsilon_i(f)) \tilde{\pi}(0, \epsilon_i(f))^* = 0$, and $ff^* =
0$ if and only if $f = 0$.  

If $\tilde{\pi}(0, \epsilon_i(f)) = 0$ and $f \geq 0$, then 
\begin{align*}
\| \tilde{t} (0, \epsilon_i(\sqrt f)) \|^2 & = \| \tilde{t} (0,
\epsilon_i(\sqrt f))^* \tilde{t} (0, \epsilon_i(\sqrt f)) \| \\ 
& = \| \tilde{\pi} ( \langle (0, \epsilon_i(\sqrt f)), (0, \epsilon_i(\sqrt f))
\rangle_B \| \\ & = \| \tilde{\pi}  (0, \epsilon_i(f)) \| \\ 
& = 0
\end{align*}
so that $\tilde{t} (0, \epsilon_i(\sqrt f)) = 0$ and consequently 
\begin{align*}
0 & = \tilde{t} (0, \epsilon_i(\sqrt f)) \tilde{t} (0, \epsilon_i(\sqrt
f))^* = \tilde{\pi}^{(1)} (\Theta_{(0, \epsilon_i(\sqrt f)), (0,
\epsilon_i(\sqrt f))}^Y) \\ & =  \begin{cases} \tilde{\pi}^{(1)}
(\phi_B(f,\vec{0}))& \text{ if $i=1$} \\ \tilde{\pi}^{(1)}(\phi_B(0,
\epsilon_{i-1}(f))) & \text{ if $i \geq 2$} \end{cases} = \begin{cases}
\tilde{\pi} (f,\vec{0}) & \text{ if $i=1$} \\ \tilde{\pi} (0,
\epsilon_{i-1}(f)) & \text{ if $i \geq 2$.} \end{cases}
\end{align*}
If $i=1$, the fact that $\tilde{\pi}|_A$ is injective implies that $f = 0$. 
Furthermore, an inductive argument combined with the above equality shows
that for all $i \in \N$ we have $f = 0$.
\end{proof}

\begin{lemma} \label{final-space}
Let $(\tilde{\pi}, \tilde{t})$ be a representation of $Y$ which is coisometric
on $J_Y = J(Y)$.  If $\vec{f} = (f_1, f_2, \ldots ) \in T$ and $\vec{g} = (g_1,
g_2, \ldots) \in T$, then $$\tilde{t}(0, \vec{f}) \tilde{t}(0, \vec{g})^* =
\tilde{\pi}(f_1g_1^*, (f_2g_2^*, f_3g_3^*, \ldots)).$$
\end{lemma}

\begin{proof}
For any $(\xi, \vec{h}) \in Y = X \oplus T$ we have
$$\Theta_{(0,\vec{f}),(0,\vec{g})} (\xi, \vec{h}) = (0,\vec{f})
\langle (0,\vec{g}), (\xi, \vec{h}) \rangle_B = (0, \vec{f} \vec{g}^* \vec{h})
= \phi_B (f_1g_1^*, ((f_2g_2^*,  \ldots)) (\xi, \vec{h})$$ so that
$\Theta_{(0,\vec{f}),(0,\vec{g})} = \phi_B (f_1g_1^*, (f_2g_2^*,  \ldots))$. 
Thus 
\begin{align*}
\tilde{t}(0, \vec{f}) \tilde{t}(0, \vec{g})^* & = \tilde{\pi}^{(1)}
(\Theta_{(0,\vec{f}),(0,\vec{g})}) = \tilde{\pi}^{(1)} (\phi_B (f_1g_1^*,
(f_2g_2^*, \ldots)) \\
& = \tilde{\pi}(f_1g_1^*, (f_2g_2^*, f_3g_3^*, \ldots)).
\end{align*} 
\end{proof}

\begin{lemma}
\label{im-relations} Let $(\tilde{\pi}, \tilde{t})$ be a representation of
$Y$.  If $\xi \in X$, $a \in A$, and $\vec{f} \in T$, then the following
relations hold: 
\begin{align*}
(1) & \ \ \tilde{t}(0,\vec{f})\tilde{\pi}(a,\vec{0}) = 0 \\
(2) & \ \ \tilde{t}(0,\vec{f})\tilde{t}(\xi,\vec{0}) = 0 \\
(3) & \ \ \tilde{t}(0,\vec{f})\tilde{t}(\xi,\vec{0})^* =0
\end{align*}
\end{lemma}

\begin{proof}
To see (1) we note that $\tilde{t}(0,\vec{f})\tilde{\pi}(a,\vec{0}) =
\tilde{t}((0,\vec{f} )(a,\vec{0})) = \tilde{t}(0,0) =0$. To see (2) and (3) let
$\{ \vec{e}_\lambda
\}_{\lambda \in \Lambda}$ be an approximate unit for $T$. Then 
\begin{equation*}
\tilde{t}(0,\vec{f}) \tilde{t}(\xi,\vec{0}) = \lim_\lambda \tilde{t}(0,
\vec{f}\vec{e}_\lambda) \tilde{t}(\xi,\vec{0}) = \lim_\lambda \tilde{t}(0,
\vec{f}) \tilde{\pi}(0,\vec{e}_\lambda) \tilde{t}(\xi,\vec{0}) = 0
\end{equation*}
which shows that (2) holds, and 
\begin{align*}
\tilde{t}(0,\vec{f}) \tilde{t}(\xi, \vec{0})^* &= \lim_\lambda \tilde{t}(0,
\vec{f}\vec{e}_\lambda) \tilde{t}(\xi,\vec{0})^* = \lim_\lambda \tilde{t}(0,
\vec{f}) \tilde{\pi}(0,\vec{e}_\lambda) \tilde{t}(\xi,\vec{0})^* \\ &=
\lim_\lambda
\tilde{t}(0,
\vec{f}) ( \tilde{t}(\xi,\vec{0}) \tilde{\pi}(0,\vec{e}_\lambda) )^* = 0
\end{align*}
which shows that (3) holds.
\end{proof}

\begin{lemma} \label{im-perp} 
Let $(\tilde{\pi}, \tilde{t})$ be a representation of $Y$, and define $(\pi,t)
:= (\tilde{\pi}|_A, \tilde{t}|_X)$.  If $c \in C^*(\pi,t)$ and $\vec{f}
\in T$, then
\begin{equation*}
\tilde{t}(0,\vec{f}) c = 0.
\end{equation*}
\end{lemma}

\begin{proof}
Since $C^*(\pi,t)$ is generated by elements of the form $\tilde{\pi}(a,
\vec{0})$ and $\tilde{t}(\xi,\vec{0})$, the result follows from the relations in
Lemma~\ref{im-relations}.
\end{proof}

\begin{lemma} \label{core-in-corner} 
Let $(\tilde{\pi}, \tilde{t})$ be a representation of $Y$ which is coisometric
on $J_Y = J(Y)$, and define
$(\pi,t) := (\tilde{\pi}|_A, \tilde{t}|_X)$.   If $n \in \{0, 1, 2, \ldots \}$,
then any element of the form
\begin{equation*}
\tilde{t}(\xi_1, \vec{f_1}) \ldots \tilde{t}(\xi_n,\vec{f_n}) \tilde{\pi}(a,
\vec{h})
\tilde{t}(\eta_n, \vec{g_n})^* \ldots \tilde{t}(\eta_1,\vec{g_1})^*
\end{equation*}
will be equal to $c + \tilde{\pi}(0,\vec{k})$ for some $c \in
C^*(\pi,t)$ and some $\vec{k} \in T$.
\end{lemma}

\begin{proof}
We shall prove this by induction on $n$.

\noindent \textsc{Base Case:} $n=0$. Then the term above is equal to
$\tilde{\pi}(a, 
\vec{h}) = \tilde{\pi}(a, \vec{0}) + \tilde{\pi}(0,\vec{h})$ and the
claim holds trivially.

\noindent \textsc{Inductive Step:} Assume the claim holds for $n$. Given an
element 
\begin{equation*}
\tilde{t}(\xi_1, \vec{f_1}) \ldots \tilde{t}(\xi_{n+1},\vec{f}_{n+1})
\tilde{\pi}(a,
\vec{h}) \tilde{t}(\eta_{n+1},\vec{g}_{n+1})^* \ldots
\tilde{t}(\eta_1,\vec{g_1})^*
\end{equation*}
it follows from the inductive hypothesis that 
\begin{equation*}
\tilde{t}(\xi_2, \vec{f_2}) \ldots \tilde{t}(\xi_{n+1},\vec{f}_{n+1})
\tilde{\pi}(a,
\vec{h}) \tilde{t}(\eta_{n+1},\vec{g}_{n+1})^* \ldots
\tilde{t}(\eta_2,\vec{g_2})^*
\end{equation*}
has the form $c + \tilde{\pi}(0,\vec{k})$ for $c \in
C^*(\pi,t)$ and $\vec{k} \in T$. Thus using Lemma~\ref{im-perp} gives 
\begin{align*}
& \tilde{t}(\xi_1, \vec{f_1}) \ldots \tilde{t}(\xi_{n+1},\vec{f}_{n+1})
\tilde{\pi}(a,
\vec{h}) \tilde{t}(\eta_{n+1},\vec{g}_{n+1})^* \ldots
\tilde{t}(\eta_1,\vec{g_1})^* \\ 
= \ &\tilde{t}(\xi_1,
\vec{f_1}) (c + \tilde{\pi}(0,\vec{k})) \tilde{t}(\eta_1,\vec{g_1} )^*
\\ = \ &(\tilde{t}(\xi_1, \vec{0}) + \tilde{t}(0, \vec{f_1})) (c +
\tilde{\pi}(0,\vec{k})) (\tilde{t}(\eta_1, \vec{0}) + \tilde{t}(0,
\vec{g_1})^*) \\ = \ &\tilde{t}(\xi_1,\vec{0}) c \tilde{t}(\eta_1,
\vec{0})^* +
\tilde{t}(0,\vec{f_1})\tilde{\pi}(0, \vec{k})\tilde{t}(0,\vec{g_1})^*
\\ = \ &\tilde{t}(\xi_1,\vec{0})
c \tilde{t}(\eta_1, \vec{0})^* + \tilde{t}(0,\vec{f_1} \vec{k})
\tilde{t}(0,\vec{g_1})^*.
\end{align*}
It follows from Lemma~\ref{final-space} that $\tilde{t}(0,\vec{f_1}
\vec{k}) \tilde{t}(0,\vec{g_1})^*$ is of the form $c^{\prime}+
\tilde{\pi}(0,\vec{k^{\prime}})$ with $c^{\prime}\in \im \pi \subseteq
C^*(\pi,t)$. Since
$\tilde{t}(\xi_1,\vec{0}) c \tilde{t}(\eta_1, \vec{0})^*$ is also in
$C^*(\pi,t)$ the proof is complete.
\end{proof}

We wish to show that if $(\tilde{\pi},\tilde{t})$ is a representation of $Y$
and if we restrict to obtain $(\pi,t) := (\tilde{\pi}|_X,\tilde{t}|_A)$, then
$C^*(\pi,t)$ is a corner of $C^*(\tilde{\pi},\tilde{t})$.  If $A$ is unital and
$X$ is left essential, then this corner will be determined by the
projection $\pi(1,\vec{0})$. However, in the following lemma we wish to consider
the general case and must make use of approximate units to define the
projection that determines the corner.

\begin{lemma} \label{proj-in-OY}
Let $X$ be a $C^*$-correspondence over $A$ and let $T := (\ker \phi)^\N$ be the
tail determined by $\ker \phi$.  If $Y := X \oplus T$ is the
$C^*$-correspondence over $B := A \oplus T$ formed by adding the tail $T$ to
$X$, and if $(\tilde{\pi}, \tilde{t})$ is a representation of $Y$, then there
exists a projection $p \in \mathcal{M} (C^*(\tilde{\pi},\tilde{t}))$ with the
property that for all $a \in A$, $\xi \in X$, and $\vec{f} \in T$ the following
relations hold: 
\begin{align*}
(1) & \ \ p\tilde{t}(\xi,\vec{f}) = \tilde{t}(\xi, (f_1, 0, 0,
\ldots)) \\ (2) & \ \ \tilde{t}(\xi,\vec{f})p = \tilde{t}(\xi, \vec{0}) \\
(3) & \ \ p\tilde{\pi}(a , \vec{f}) = \tilde{\pi}(a , \vec{f})p =
\tilde{\pi}(a, \vec{0})
\end{align*}
\end{lemma}

\begin{proof}
Let $\{ \vec{e}_\lambda \}_{\lambda \in \Lambda}$ be an approximate unit for
$T$, and for each $\lambda \in \Lambda$ let $\vec{e}_\lambda = (e_\lambda^1,
e_\lambda^2, \ldots)$. Consider $\{\tilde{\pi}(0, \vec{e}_\lambda) \}_{\lambda
\in \Lambda}$. For any element 
\begin{equation}  \label{typ-element}
\tilde{t}(\xi_1, \vec{f_1}) \ldots \tilde{t}(\xi_n,\vec{f_n}) \tilde{\pi}(a,
\vec{h}) \tilde{t}(\eta_m, \vec{g_m})^* \ldots \tilde{t}(\eta_1,\vec{g_1})^*
\end{equation}
we have 
\begin{align*}
&\lim_\lambda \tilde{\pi}(0,\vec{e}_\lambda) \tilde{t}(\xi_1, \vec{f_1}) \ldots
\tilde{t}(\xi_n, \vec{f_n}) \tilde{\pi}(a, \vec{h})
\tilde{t}(\eta_m,\vec{g_m})^*
\ldots \tilde{t}(\eta_1,\vec{ g_1})^* \\ 
= \ & \lim_\lambda \tilde{t}( 0, (0, e_\lambda^1 f_{12}, e_\lambda^2 f_{13},
\ldots) \ldots \tilde{t}(\xi_n,\vec{f_n}) \tilde{\pi}(a, \vec{h})
\tilde{t}(\eta_m,\vec{g_m})^* \ldots \tilde{t}(\eta_1,\vec{g_1})^* \\ 
= \ & \tilde{t}(0, (0, f_{12}, f_{13}, \ldots)) \ldots
\tilde{t}(\xi_n,\vec{f_n}) \tilde{\pi}(a, \vec{h})
\tilde{t}(\eta_m,\vec{g_m})^* \ldots \tilde{t}(\eta_1,\vec{g_1})^*
\end{align*}
so this limit exists.

Now since any $c \in C^*(\tilde{\pi},\tilde{t})$ can be approximated by a
finite sum of elements of the form shown in (\ref{typ-element}), it follows that
$\lim_\lambda \tilde{\pi}(0, \vec{e}_\lambda) c$ exists for all $c
\in C^*(\tilde{\pi}, \tilde{t})$. Let us view $C^*(\tilde{\pi}, \tilde{t})$ as
a $C^*$-correspondence over itself (see \cite[Example~2.10]{RW}). If we define
$q : C^*(\tilde{\pi}, \tilde{t}) \to  C^*(\tilde{\pi}, \tilde{t})$ by $q(c)
= \lim_\lambda \tilde{\pi}(0, \vec{e}_\lambda) c$ then we see that for any
$c, d \in C^*(\tilde{\pi}, \tilde{t})$ we have 
\begin{equation*}
d^* q (c) = \lim_\lambda d^* \tilde{\pi}(0, \vec{e}_\lambda)
c = \lim_\lambda ( \tilde{\pi}(0, \vec{e}_\lambda) d)^* c =
q(d)^* c
\end{equation*}
and hence $q$ is an adjointable operator on $C^*(\tilde{\pi}, \tilde{t})$.
Therefore $q$ defines (left multiplication by) an element
in the multiplier algebra $\mathcal{M}(C^*(\tilde{\pi}, \tilde{t}))$
\cite[Theorem~2.47]{RW}. It is easy to check that $q^2 = q^*=q$ so that $q$ is a
projection. Now if we let $p:= 1-q$ in $\mathcal{M} (C^*(\tilde{\pi},
\tilde{t}))$, then it is easy to check that relations (1), (2), and (3) follow
from the definition of $q$.
\end{proof}

\noindent \emph{Proof of Theorem~\ref{pi-t-extend}.} (a)  Let $I:= \ker
\phi$, set $\Hi_0 := \pi(I) \Hi_X$, and define $\Hi_T :=
\bigoplus_{i=1}^\infty \Hi_i$ where $\Hi_i = \Hi_0$ for all $i = 1,2,
\ldots$.  We define $\tilde{t} : Y \to \B (\Hi_X \oplus \Hi_T)$ and
$\tilde{\pi} : B \to \B (\Hi_X \oplus \Hi_T)$ as follows:  Viewing $Y$ as $Y
= X \oplus T$ and $B$ as $B = A \oplus T$, for any $(h,(h_1,h_2, \ldots ))
\in \Hi_\Q \oplus \Hi_T$ we define $$\tilde{t}(\xi, (f_1, f_2,
\ldots))(h,(h_1,h_2, \ldots )) = (t(\xi)h + \pi(f_1)h_1, (\pi(f_2)h_2,
\pi(f_3)h_3, \ldots))$$ and $$\tilde{\pi} (a, (f_1, f_2, \ldots))(h,(h_1,h_2,
\ldots )) = (\pi(a)h,(\pi(f_1)h_1,\pi(f_2)h_2, \ldots )).$$  Then it is
straightforward to show that $(\tilde{\pi}, \tilde{t})$ is a representation of
$Y$ on $\Hi_\Q \oplus \Hi_T$.  To see that $(\tilde{\pi},t)$ is coisometric on
$J(Y)$, choose an element $(a, (f_1, f_2, \ldots)) \in J(Y)$.  By
Lemma~\ref{identify-J(Y)} we know that $a = a_1 + a_2$ for $a_1 \in J_X$ and
$a_2 \in \ker \phi$.  Furthermore, since $a_1 \in J(X)$ we may write
$\phi(a_1) = \lim_n \sum_{k=1}^{N_n}
\Theta^X_{\xi_{n,k}, \eta_{n,k}}$ for some $\xi_{n,k}, \eta_{n,k} \in X$. 
It follows that $$\phi_B(a_1, \vec{0}) = \lim_n \sum_{k=1}^{N_n}
\Theta^Y_{(\xi_{n,k}, \vec{0}), (\eta_{n,k}, \vec{0})} \in \K(Y).$$ In
addition, since $a_2 \in \ker \phi$ we see that if we let $\{ \vec{e}_\lambda
\}_{\lambda \in \Lambda}$ be an approximate unit for $T$ with $\vec{e}_\lambda =
(e_\lambda^1, e_\lambda^2, \ldots )$ for each $\lambda$, then $$\phi_B(a_2
,\vec{f}) = \lim_\lambda \Theta^Y_{(0, (a, f_1, f_2, \ldots)), (0,
(e_\lambda^1, e_\lambda^2, \ldots))} \in \K(Y).$$  Now for any $n \in \N$ we
see that $\{e_\lambda^n \}_{\lambda \in \Lambda}$ is an approximate unit for
$\ker \phi$.  Furthermore, we see that for all $(\xi,
\vec{f}), (\eta, \vec{g}) \in Y = X \oplus T$ we have $$\tilde{t}(\xi,
\vec{f})\tilde{t}(\eta, \vec{g})^* = (t(\xi)t(\eta)^* +
\pi(f_1g_1^*), (\pi(f_2g_2^*), \pi(f_3g_3^*), \ldots ))$$ and thus 
\begin{align*}
& \ \tilde{\pi}^{(1)} (\phi_B(a, \vec{f})) \\
= & \ \tilde{\pi}^{(1)}(\phi_B(a_1, \vec{0})) + \tilde{\pi}^{(1)}(\phi_B(a_2, \vec{f})) \\
= & \ \lim_n \sum_{k=1}^{N_n} \tilde{t}(\xi_{n,k}, \vec{0}) \tilde{t}
(\eta_{n,k}, \vec{0})^* + \lim_\lambda \tilde{t} (0, (a_2, f_1, f_2,
\ldots)) \tilde{t} (0, (e_\lambda^1, e_\lambda^2, \ldots))^* \\ 
= & \ \lim_n \sum_{k=1}^{N_n} (t(\xi_{n,k}) t(\eta_{n,k})^*, \vec{0}) +
\lim_\lambda (\pi(a_2 e_\lambda^1), (\pi(f_1 e_\lambda^2), \pi(f_2
e_\lambda^3), \ldots)) \\ = & \ (\pi^{(1)} (\phi(a_1)), \vec{0}) + (\pi(a_2),
(\pi(f_1), \pi(f_2), \ldots)) \\ 
= & \ (\pi(a_1), \vec{0}) + (\pi(a_2), (\pi(f_1), \pi(f_2), \ldots)) \\ = & \
\tilde{\pi}(a, \vec{f}) 
\end{align*}
so $(\tilde{\pi}, \tilde{t})$ is coisometric on $J(Y)$. \\

\noindent (b)  If $(\tilde{\pi}, \tilde{t})$ is a representation of
$Y$ in a $C^*$-algebra $C$ which is coisometric on $J(Y)$, then it is
straightforward to see that the restriction $(\tilde{\pi}|_A, \tilde{t}|_X)$
is a representation.  To see that $(\tilde{\pi}|_A, \tilde{t}|_X)$
is coisometric on $J_X$, choose an element $a \in J_X$.  Since $J_X \subseteq
J(X)$ we may write $\phi(a) = \lim_n \sum_{k=1}^{N_n} \Theta^X_{\xi_{n,k},
\eta_{n,k}}$ for some $\xi_{n,k}, \eta_{n,k} \in X$.  In addition, since $a \in
(\ker \phi)^\perp \subseteq J_X$ we have that
$$\phi_B(a,
\vec{0}) = \lim_n \sum_{k=1}^{N_n} \Theta^Y_{(\xi_{n,k}, \vec{0}), (\eta_{n,k},
\vec{0})} \in \K(Y)$$ and we have
\begin{align*}
\tilde{\pi}|_A^{(1)} (\phi_A(a)) & = \lim_n \sum_{k=1}^{N_n}
\tilde{t}|_X(\xi_{n,k}) \tilde{t}|_X(\eta_{n,k})^*  =  \lim_n
\sum_{k=1}^{N_n} \tilde{t}(\xi_{n,k}, \vec{0}) \tilde{t}(\eta_{n,k},
\vec{0})^* \\ 
& = \tilde{\pi}^{(1)}(\phi_B(a, \vec{0})) = \tilde{\pi} (a,
\vec{0}) = \tilde{\pi}|_A(a)
\end{align*}
so $(\tilde{\pi}|_A, \tilde{t}|_X)$ is coisometric on $J_X$.

Furthermore, suppose that the restriction $\tilde{\pi}|_A$ is injective.  If
$(a, \vec{f}) \in B := A \oplus T$ and $\tilde{\pi}(a,\vec{f}) = 0$, let
$\{g_\lambda \}_{\lambda \in \Lambda}$ be an approximate unit for
$\ker \phi$, and for any $f \in \ker \phi$ and $i \in \N$ let $\epsilon_i (f) :=
(0, \ldots, 0, f, 0, \ldots)$ where $f$ is in the $i$\textsuperscript{th}
position.  Since $\tilde{\pi}(a,\vec{f}) = 0$ we see that if we write $\vec{f}
= (f_1, f_2, \ldots)$, then for all $i \in \N$ we have
$$\tilde{\pi}(0,\epsilon_i(g_\lambda f_i)) = \tilde{\pi} (0,
\epsilon_i(g_\lambda)) \tilde{\pi}(a, \vec{f}) = 0,$$ and taking limits with
respect to $\lambda$ shows that $\tilde{\pi}(0, \epsilon_i(f_i)) =0$ for all $i
\in \N$.  From Lemma~\ref{b-inj-lemma} it follows that $f_i = 0$ for all $i \in
\N$.  Thus $\vec{f} = 0$, and since $\tilde{\pi}|_A$ is injective we also have
that $a = 0$.  Hence $\tilde{\pi}$ is injective. \\

\noindent (c)  The fact that $(\pi, t) := (\pi_Y|_A, t_Y|_X)$ is
a representation which is coisometric on $J_X$ follows from Part~(b). 
Furthermore, the fact that $\rho_{(\pi,t)}$ is injective follows from
Part~(a) which shows that any $*$-representation of $\mathcal{O}_X$ factors
through a $*$-representation of $\mathcal{O}_Y$.  All that remains is to show
that $\operatorname{im} \rho_{(\pi, t)} = C^*(\pi,t)$ is a full corner of
$\mathcal{O}_Y$.

Let $p \in \mathcal{M}(\mathcal{O}_Y)$ be the projection described in
Lemma~\ref{proj-in-OY}. We shall first show that $C^*(\pi,t) = p 
\mathcal{O}_Y p$. To begin, we see from the relations in
Lemma~\ref{proj-in-OY} that for all $a \in A$ we have $p \pi(a) p = p
\pi_Y(a, \vec{0}) p = \pi_Y (a,\vec{0}) = \pi(a)$ and for all $\xi \in X$ we
have $p t(\xi) p = p (t_Y (\xi, \vec{0}))p = t_Y(\xi, \vec{0})
= t(\xi)$. Thus $C^*(\pi,t) \subseteq p\mathcal{O}_Y p$.

To see the reverse inclusion, note that any element in $\mathcal{O}_Y$ is
the limit of sums of elements of the form 
\begin{equation*}
t_Y(\xi_1, \vec{f_1}) \ldots t_Y(\xi_n,\vec{f_n}) \pi_Y(a, \vec{h})
t_Y(\eta_m, \vec{g_m})^* \ldots t_Y(\eta_1,\vec{g_1})^*
\end{equation*}
and thus any element of $p\mathcal{O}_Y p$ is the limit of sums of elements of
the form 
\begin{equation*}
p t_Y(\xi_1, \vec{f_1}) \ldots t_Y(\xi_n,\vec{f_n}) \pi_Y(a, \vec{h})
t_Y(\eta_m,\vec{g_m})^* \ldots t_Y(\eta_1,\vec{g_1})^*p
\end{equation*}
Therefore, it suffices to show that each of these elements is in $C^*(\pi,t)$.
Now if $n \geq m$, then we may use Lemma~\ref{core-in-corner} to write 
\begin{equation*}
t_Y(\xi_{n-m+1}, \vec{f}_{n-m+1}) \ldots t_Y(\xi_n,\vec{f_n}) \pi_Y(a,
\vec{h}) t_Y(\eta_m,\vec{g_m})^* \ldots t_Y(\eta_1,\vec{g_1})^*
\end{equation*}
as $c + \pi_Y(0,\vec{k})$ for $c \in C^*(\pi,t)$ and $\vec{k} \in T$.
Then 
\begin{align*}
& p t_Y(\xi_1, \vec{f_1}) \ldots t_Y(\xi_n,\vec{f_n}) \pi_Y(a, \vec{h})
t_Y(\eta_m,\vec{g_m})^* \ldots t_Y(\eta_1,\vec{g_1})^*p \\
= \ & p t_Y(\xi_1, \vec{f_1}) \ldots t_Y(\xi_{n-m},\vec{f}_{n-m}) (c +
\pi_Y(0,\vec{k}) ) p \\
= \ & p t_Y(\xi_1, \vec{f_1}) \ldots t_Y(\xi_{n-m},\vec{f}_{n-m}) c p \\
= \ & p t_Y(\xi_1, \vec{f_1}) \ldots t_Y(\xi_{n-m},\vec{f}_{n-m}) p c p
\\
= \ & p t_Y(\xi_1, \vec{f_1}) \ldots t_Y(\xi_{n-m-1},\vec{f}_{n-m-1}) p
t_Y(\xi_{n-m},\vec{0}) p c p \\
\ \vdots & \\
= \ & p t_Y(\xi_1, \vec{0})p \ldots pt_Y(\xi_{n-m},\vec{0}) p c p \\
= \ & t_Y(\xi_1, \vec{0}) \ldots t_Y(\xi_{n-m},\vec{0}) c \\
= \ & t(\xi_1) \ldots t(\xi_{n-m}) c
\end{align*}
which is in $C^*(\pi,t)$. The case when $n \leq m$ is
similar. Hence $p \mathcal{O}_Y p \subseteq C^*(\pi,t)$.

To see that the corner $C^*(\pi,t) = p \mathcal{O}_Y p$ is
full, suppose that $\mathcal{I}$ is an ideal in $\mathcal{O}_Y$ that
contains $C^*(\pi,t)$. For $f \in \ker \phi$ and $n \in N$
define $\epsilon_n(f) := (0, \ldots, 0, f, 0, \ldots) \in T$,
where the term $f$ is in the $n^{\text{th}}$ position. Let $\{
e_\lambda \}_{\lambda \in \Lambda}$ be an approximate unit for $\ker \phi$. Now
$t_Y(\xi,\vec{0}), \pi_Y(a,\vec{0}) \in C^*(\pi,t) \subseteq 
\mathcal{I}$ for all $a \in A$ and $\xi \in X$, and since $T$ is the
$c_0$-direct sum of countably many copies of $\ker \phi$ in order to show that
$\mathcal{I}$ is all of $\mathcal{O}_Y$ it suffices to prove that for all $n
\in N$ and $\lambda \in
\Lambda$ we have $t_Y(0, \epsilon_n(e_\lambda)) \in 
\mathcal{I}$ and $\pi_Y (0, \epsilon_n(e_\lambda)) \in \mathcal{I}$. We shall
prove this by induction on $n$.

\noindent \textsc{Base Case:} For any $\beta, \lambda \in \Lambda$ we have from
Lemma~\ref{final-space} that
\begin{equation*}
t_Y(0, \epsilon_1(e_\lambda)) t_Y(0,\epsilon_1(e_\beta))^* =
\pi_Y(e_\lambda e_\beta, \vec{0}) \in \mathcal{I}.
\end{equation*}
Also for any $\alpha \in \Lambda$ we have 
\begin{align*}
t_Y(0, (\epsilon_1(e_\lambda e_\beta e_\alpha)) &= t_Y(0,
\epsilon_1(e_\lambda)) \pi_Y(0,\epsilon_1(e_\beta^* e_\alpha)) \\
&= t_Y(0, \epsilon_1(e_\lambda)) t_Y(0,\epsilon_1(e_\beta))^* t_Y(0,
\epsilon_1(e_\alpha))
\end{align*}
which is in $\mathcal{I}$. Taking limits with respect to $\alpha$ and $\beta$
gives 
\begin{equation*}
t_Y(0,\epsilon_1(e_\lambda)) = \lim_\beta \lim_\alpha t_Y(0,
\epsilon_1(e_\lambda e_\beta e_\alpha)) \in \mathcal{I}.
\end{equation*}
Furthermore, since $t_Y(0,\epsilon_1(e_\lambda)) \in \mathcal{I}$ for all
$\lambda \in
\Lambda$, we see that 
\begin{equation*}
\pi_Y(0,\epsilon_1(e_\lambda)) = \lim_\beta \pi_Y(0, \epsilon_1(e_\lambda
e_\beta)) = \lim_\beta t_Y(0,\epsilon_1(e_\lambda))^*
t_Y(0,\epsilon_1(e_\beta)) \in \mathcal{I}.
\end{equation*}

\noindent \textsc{Inductive step:} Suppose that $t_Y(0,\epsilon_n(e_\lambda)),
\pi_Y(0,\epsilon_n(e_\lambda)) \in \mathcal{I}$ for any $\lambda
\in \Lambda$. Then for all $\lambda, \beta \in \Lambda$ we have 
\begin{align*}
t_Y(0,\epsilon_{n+1}(e_\lambda)) t_Y(0,\epsilon_{n+1}(e_\beta))^* &=
\pi_Y^{(1)} (\Theta_{(0,\epsilon_{n+1}(e_\lambda)),
(0,\epsilon_{n+1}(e_\beta))}^Y) \\
&= \pi_Y^{(1)}(\phi_B(0, \epsilon_n(e_\beta e_\lambda)) ) \\
&= \pi_Y(0,\epsilon_n(e_\beta e_\lambda)) \\
&= \pi_Y(0, \epsilon_n(e_\beta)) \pi_Y(0, \epsilon_n(e_\lambda))
\end{align*}
which is in $\mathcal{I}$. Thus for any $\alpha \in \Lambda$ we have that 
\begin{align*}
t_Y(0,\epsilon_{n+1}(e_\lambda e_\beta e_\alpha)) &=
t_Y(0,\epsilon_{n+1}(e_\lambda)) \pi_Y(0,\epsilon_{n+1}(e_\beta e_\alpha)) \\
&= t_Y(o,\epsilon_{n+1}(e_\lambda)) t_Y(0,\epsilon_{n+1}(e_\beta))^*
t_Y(0,\epsilon_{n+1}(e_\alpha))
\end{align*}
is in $\mathcal{I}$. Taking limits with respect to $\alpha$ and $\beta$
gives 
\begin{equation*}
t_Y(0,\epsilon_{n+1}(e_\lambda)) = \lim_\beta \lim_\alpha t_Y(0,
\epsilon_{n+1}(e_\lambda e_\beta e_\alpha)) \in \mathcal{I}.
\end{equation*}
Furthermore, since $t_Y(0,\epsilon_{n+1}(e_\lambda)) \in \mathcal{I}$ for
all $\lambda \in \Lambda$, we have 
\begin{equation*}
\pi_Y(0, \epsilon_{n+1}(e_\lambda)) = \lim_\beta \pi_Y(0, \epsilon_{n+1}(e_\beta
e_\lambda )) = \lim_\beta t_Y(0, \epsilon_{n+1}(e_\beta))^* t_Y(0,
\epsilon_{n+1}(e_\lambda)) \in \mathcal{I}.
\end{equation*}
\hfill \qed

%%%%%%%%%%%%%%%%%%%%%%%%%%%%%%%%%%%%%%%%%%%%%%%%%%%%%%%%%%
\section{Gauge-Invariant Uniqueness} \label{GIU-sec}
%%%%%%%%%%%%%%%%%%%%%%%%%%%%%%%%%%%%%%%%%%%%%%%%%%%%%%%%%%

Recall that we let $\gamma$ denote the gauge action of $\T$ on
$\mathcal{O}_X$.   A gauge-invariant uniqueness was proven in
\cite[Theorem~4.1]{FMR} for (augmented) Cuntz-Pimsner algebras.  Our method of
adding tails, together with Theorem~\ref{pi-t-extend}, will allow us to extend
this theorem to the case when $\phi$ is not injective, and ultimately to all
relative Cuntz-Pimsner algebras.  

The following Gauge-Invariant Uniqueness Theorem was proven by Katsura using
direct methods in \cite[Theorem~6.4]{Kat4}.  We shall now give an alternate
proof, showing how the method of adding tails can be used to bootstrap
\cite[Theorem~4.1]{FMR} to the general case.

\begin{theorem}[Gauge-Invariant Uniqueness] \label{GIU} Let $X$ be a
$C^*$-correspondence over $A$, and let $(\pi_X, t_X)$ be a universal $J(X)$-coisometric
representation of $X$. If $\rho : \mathcal{O}_X
\rightarrow C$ is a homomorphism between $C^*$-algebras which satisfies the
following two conditions:
\begin{enumerate}
\item the restriction of $\rho$ to $\pi_X(A)$ is injective
\item there is a strongly continuous action $\beta : \mathbb{T}
\rightarrow \operatorname{Aut}(\rho(\mathcal{O}_X))$ such that $\beta_z \circ
\rho = \rho \circ \gamma_z$ for all $z \in \mathbb{T}$
\end{enumerate}
then $\rho$ is injective.
\end{theorem}

\begin{remark}
When $\phi$ is injective, the statement above is actually an
equivalent reformulation of \cite[Theorem~4.1]{FMR}.  The equivalence relies on
the fact that for any $C^*$-correspondence $X$, the universal $J(X)$-coisometric representation
$(i_A, i_X)$ has the property that $i_A$ is
injective if and only if the left action $\phi$ is injective.
\end{remark}

\noindent \emph{Proof of Theorem~\ref{GIU}.}  Let $T := (\ker \phi)^\N$ be the
tail determined by $\ker \phi$, and let $Y := X \oplus T$ be the
$C^*$-correspondence over $B := A \oplus T$ formed by adding the tail $T$ to
$X$.  By Theorem~\ref{pi-t-extend}(c) we may identify $(\mathcal{O}_X, \pi_X,
t_X)$ with
$(S, \pi_Y|_A, t_Y|_X)$ where $S$ is the $C^*$-subalgebra of $\mathcal{O}_Y$
generated by $$\{ \pi_Y( a, \vec{0}), t_Y(\xi, \vec{0}) : a \in A \text{ and } 
\xi \in X \}.$$

Since $\beta : \T  \to \aut(\im \rho)$ is an action of $\T$ on $\im \rho$, there
exists a Hilbert space $\Hi_X$, a faithful representation $\kappa : \im \rho \to
\B (\Hi_X)$, and a unitary representation $U : \T \to \mathcal{U}(\Hi_X)$ such
that $$\kappa (\beta_z(x)) = U_z \kappa(x) U_z^* \qquad \text{ for all $x \in
\im \rho$ and $z \in \T$.}$$  In addition, since $\tau := \kappa \circ \rho$ is
a $*$-homomorphism from $S$ into $\B (\Hi_X)$ which is faithful on
$\pi_X(A)$, it follows from Theorem~\ref{pi-t-extend}(a) that $\tau$ may be
extended to a $*$-homomorphism $\tilde{\tau} : \mathcal{O}_Y \to \B (\Hi_X
\oplus \Hi_T)$ with $\tilde{\tau}$ faithful on $\pi_Y (B)$.

We shall now define a unitary representation $W : \T \to \B (\Hi_X \oplus
\Hi_T)$ as follows.  We see from the proof of Theorem~\ref{pi-t-extend}(a)
that $\Hi_T := \bigoplus_{i=1}^\infty \Hi_i$.  Thus for $(h, (h_1, h_2,
\ldots)) \in \Hi_\Q \oplus \Hi_T$ we define $$W_z (h, (h_1, h_2, \ldots)) :=
(U_zh, (z^{-1}h_1, z^{-2}h_2, \ldots)) \qquad \text{ for $z \in \T$.}$$  We may
then define $\tilde{\beta} : \T \to \aut (\B  (\Hi_X \oplus \Hi_T))$ by
$\tilde{\beta}_z (T_0) := W_z T_0 W_z^*$, and we see that $\tilde{\beta}$ is a
strongly continuous gauge action.  Furthermore, if $\gamma'$
denotes the gauge action of $\T$ on $\mathcal{O}_Y$, then $\tilde{\beta}_z
\circ \tilde{\tau} = \tilde{\tau} \circ \gamma'_z$ (to see this recall how the
extension $\tilde{\tau}$ is defined in the proof of Theorem~\ref{pi-t-extend}(a)
and then simply check on the generators $\{ t_Y(\xi, \vec{f}), \pi_Y(a,\vec{g})
: \xi \in X, a \in A, \text{ and } \vec{f},\vec{g} \in T \}$).  Thus by
\cite[Theorem~4.1]{FMR} we have that  $\tilde{\tau}$ is injective.  Hence
$\tilde{\tau}|_S = \tau = \kappa \circ \rho$ is injective, and $\rho$ is
injective.
\hfil \qed

$\text{ }$

To conclude this section we shall interpret our result in the relative
Cuntz-Pimsner setting.

\begin{remark} \label{rel-GIU}
Katsura has shown in \cite{Kat5} that if $\mathcal{O}(K,X)$
is a relative Cuntz-Pimsner algebra, then there exists a $C^*$-correspondence
$X'$ with the property that $\mathcal{O}_{X'}$ is naturally isomorphic to
$\mathcal{O}(K,X)$.  Using this analysis one can obtain the following
interpretation of Theorem~\ref{GIU} for relative Cuntz-Pimsner algebras.

$\text{ }$

\noindent \textbf{Interpretation of Theorem~\ref{GIU} for Relative Cuntz-Pimsner
Algebras:} \emph{Let $X$ be a $C^*$-correspondence with
left action $\phi : X \to \Li (X)$, let $K$ be an ideal in $J(X) :=
\phi^{-1}(\K(X))$, and let $(\pi_X, t_X)$ be a universal $K$-coisometric representation of $X$.  If $\rho : \mathcal{O}_X \rightarrow C$ is a
homomorphism between $C^*$-algebras which satisfies the following three
conditions:
\begin{enumerate}
\item[(1)] the restriction of $\rho$ to $\pi_X(A)$ is injective
\item[(2)] if $\rho(\pi_X(a)) \in \rho(\pi_X^{(1)}(\K(X)))$, then $\pi_X(a) \in
\pi_X(K)$
\item[(3)] there is a strongly continuous action $\beta : \mathbb{T}
\rightarrow \operatorname{Aut}(\rho(\mathcal{O}_X))$ such that $\beta_z \circ
\rho = \rho \circ \gamma_z$ for all $z \in \mathbb{T}$
\end{enumerate}
then $\rho$ is injective.}

$\text{ }$

Finally, we mention that if we define a map $T_K : J(X) \rightarrow
\mathcal{O}(K,X)$ by 
\begin{equation*}
T_K (a) := \pi_X(a) - \pi_X^{(1)}(\phi(a))
\end{equation*}
then the equation 
\begin{align*}
T_K(a)T_K(b) & = (\pi_X(a) - \pi_X^{(1)}(\phi(a))) (\pi_X(b) -
\pi_X^{(1)}(\phi(b))) \\
& = \pi_X(ab) - \pi_X^{(1)}(\phi(a))\pi_X(b) - \pi_X(a)\pi_X^{(1)}(\phi(b)) +
\pi_X^{(1)}(\phi(ab)) \\
& = \pi_X(ab) - \pi_X^{(1)}(\phi(ab)) \\
& = T_K(ab)
\end{align*}
shows that this map is a homomorphism.  If $\pi_X$ is injective (which by
\cite[Proposition~2.21]{MS} occurs if and only if $K \cap \ker \phi =
\emptyset$), then we may replace Condition~(2) in the above statement by the
condition
\emph{
\begin{enumerate}
\item[($2'$)] the restriction of $\rho$ to $T_K(J(X))$ is injective.
\end{enumerate}}
\end{remark}

%%%%%%%%%%%%%%%%%%%%%%%%%%%%%%%%%%%%%%%%%%%%%%%%%%%%%%%%%%
\section{Gauge-Invariant Ideals} \label{ideals-sec}
%%%%%%%%%%%%%%%%%%%%%%%%%%%%%%%%%%%%%%%%%%%%%%%%%%%%%%%%%%

In this section we use Theorem~\ref{GIU} to characterize the
gauge-invariant ideals in $C^*$-algebras associated to certain correspondences. 

\begin{definition}
\label{X-inv-sat} Let $X$ be a $C^*$-correspondence over $A$. We say that an
ideal $I \triangleleft A$ is \emph{$X$-invariant} if $\phi (I) X \subseteq XI
$. We say that an $X$-invariant ideal $I \triangleleft A$ is \emph{$X$-saturated} if 
\begin{equation*}
a \in J_X \text{ and } \phi(a)X \subseteq XI \ \Longrightarrow \ a \in I.
\end{equation*}
\end{definition}

\begin{remark}
In \cite{KPW} the authors only considered Hilbert bimodules (i.e.
$C^*$-correspondences) for which $\phi$ is injective and $\phi(A) \subseteq
\mathcal{K}(X)$, and thus the definition of $X$-saturated that they gave
was that $a \in A$ and $\phi(a)X \subseteq XI$ implies $a \in I$. Since $J_X =
A$ throughout their paper, this notion is equivalent to the one defined in
Definition~\ref{X-inv-sat}. In \cite[Remark~3.11]{FMR} it was suggested that the
definition of $X$-saturated for general $C^*$-correspondences should also be
that $a \in A$ and $\phi(a)X \subseteq XI$ implies $a \in I$. However, after
considering how the definition of \emph{saturated} was extended to (or rather
modified for) non-row-finite graphs in \cite[\S3]{BHRS} and \cite[\S3]{DT1} we
believe that Definition~\ref{X-inv-sat} is the appropriate generalization.
\end{remark}

Recall that if $I$ is an ideal of $A$, then 
\begin{equation*}
X_I := \{ x \in X : \langle x, y \rangle_A \in I \text{ for all } y \in X \}
\end{equation*}
is a right Hilbert $A$-module, and by the Hewitt-Cohen Factorization Theorem 
$X_I = XI := \{ x \cdot i : x \in X \text{ and } i \in I \}$ (see \cite[\S2]{FMR}).
Furthermore, $X / XI$ is a right Hilbert $A / I$-module in the obvious way
\cite[Lemma~2.1]{FMR}. In order for $X / XI$ to be a $C^*$-correspondence, we
need the ideal $I$ to be $X$-invariant.  Let $q^I : A \rightarrow A / I$ and
$q^{XI} : X \rightarrow X / XI$ be the appropriate quotient maps. If $I$ is
$X$-invariant, then one may define $\phi_{A/I}:A/I
\rightarrow \mathcal{L}(X/XI)$ by 
\begin{equation*}
\phi_{A/I} (q^I(a))(q^{XI}(x)) := q^{XI}(\phi(a)(x))
\end{equation*}
and with this action $X/XI$ is a $C^*$-correspondence over $A/I$
\cite[Lemma~3.2]{FMR}.

\begin{lemma} \label{row-finite-JX}
Let $X$ be a $C^*$-correspondence over a $C^*$-algebra $A$, and let $I$ be an
$X$-saturated $X$-invariant ideal in $A$.  If $q^I : A \to A / I$ denotes the
quotient map, then $$q^I(J_X) \subseteq J_{X / XI}.$$  Furthermore, if $X$ has
the following two properties:
\begin{enumerate}
\item $\phi(A) \subseteq \K (X)$
\item $\ker \phi$ is complemented in $A$ (i.e.~there exists an ideal $J$ of
$A$ with the property that $A = J \oplus \ker \phi)$,
\end{enumerate}
then $$q^I(J_X) = J_{X / XI}.$$ 
\end{lemma}

\begin{proof}
Let $a \in J_X$.  Then $a \in J(X)$, and it follows from \cite[Lemma~2.7]{FMR}
that $q^I(a) \in J(X /XI)$.  Also, if $q^I(b) \in \ker \phi_{A / I}$,
then $q^I(ab) \in \ker \phi_{A / I}$ and for all $x \in X$ we have $$q^{XI}(
\phi (ab) (x) ) = \phi_{A / I} (ab) q^{XI}(x) = 0$$ and thus 
\begin{equation} \label{sat-disp}
\phi (ab) XI \subseteq XI.
\end{equation} 
Since $a \in J_X$ and $J_X$ is an ideal, we see that $ab \in J_X$.  Now since
$I$ is $X$-saturated, (\ref{sat-disp}) implies that $ab \in I$ and
$q^I(a)q^I(b) = q^I(ab) = 0$.  Thus $q^I(a) \in (\ker \phi_{A / I})^\perp$ and
$q^I(a)
\in J_{X / XI}$.

Now suppose that Conditions (1) and (2) in the statement of the lemma hold. 
Since $\phi(A) \subseteq \K(X)$ it follows that $J(X) = A$.  In addition,
\cite[Lemma~2.7]{FMR} shows that $q^I(J(X)) = J(X / XI)$.  From Condition (2)
we know that $A = J \oplus \ker \phi$ for some ideal $J$ of $A$.  However, the
definition of $J_X$ then implies that $J=J_X$.  Thus if $a \in A$ and $q^I(a)
\in J_{X/XI}$, then we may write $a = b + c$ for $b \in J_X$ and $c \in \ker
\phi$.  But then $q^I(b) \in J_{X / XI}$ by the first part of the lemma, and
$q^I(c) = q^I(a) - q^I(b) \in J_{X / XI}$.  Since $c \in \ker \phi$ it follows
that for all $x \in X$ we have $$\phi_{A / I}(q^I(c)) q^{XI}(x) = q^{XI} (
\phi(c) (x)) = 0$$ and thus $q^I(c) \in \ker \phi_{A / I}$.  Thus $q^I(c) \in
J_{X / XI} \cap \ker \phi_{A / I} = \{ 0 \}$ so $q^I(c) = 0$ and $q^I(a) =
q^I(b) \in q^I(J_X)$.  Thus $J_{X /XI} \subseteq q^I(J_X) $.
\end{proof}

The following theorem was proven in \cite[Theorem~4.3]{KPW} under the
hypotheses that $\phi$ is injective, $A$ is unital, and $X$ is full and
finite projective as a right $A$-module (so in particular, $\phi(A)
\subseteq \mathcal{K}(X)$).  However, Theorem~\ref{GIU} allows us to
give a fairly simple proof of the result for much more general
$C^*$-correspondences. 

\begin{theorem} \label{g-i-in-OX} 
Let $X$ be a $C^*$-correspondence with the following two properties:
\begin{enumerate}
\item $\phi(A) \subseteq \K (X)$
\item $\ker \phi$ is complemented in $A$ (i.e.~there exists an ideal $J$ of
$A$ with the property that $A = J \oplus \ker \phi)$,
\end{enumerate}
and let $(\pi_X, t_X)$ be a universal $J(X)$-coisometric representation of $X$.  Then there is a lattice isomorphism from the
$X$-saturated $X$-invariant ideals of $A$ onto the gauge-invariant ideals of
$\mathcal{O}_X$ given by 
\begin{equation*}
I \mapsto \mathcal{I}(I) := \text{the ideal in $\mathcal{O}_X$ generated
by $\pi_X(I)$}
\end{equation*}
\end{theorem}

\begin{proof}
To begin we see that $\mathcal{I}(I)$ is in fact gauge invariant since 
\begin{align*}
\mathcal{I}(I) = \operatorname{\overline{\mathrm{span}}} \{ t_X(x_1) \ldots
& t_X(x_n) \pi_X(a) t_X(y_1)^* \ldots t_X(y_m)^*  \\ &: a \in I, x_1 \ldots x_n
\in X, y_1 \ldots y_m \in X, \text{and } n,m \geq 0 \}.
\end{align*}
In addition, the map $I \mapsto \mathcal{I}(I)$ is certainly inclusion
preserving.

To see that the map is surjective, let $\mathcal{I}$ be a gauge-invariant
ideal in $\mathcal{O}_X$. If we define $I := \pi_X^{-1} (\mathcal{I})$, then
it is straightforward to show that $I$ is $X$-invariant and $X$-saturated.
Now clearly $\mathcal{I}(I) \subseteq \mathcal{I}$ so there exists a
quotient map $q : \mathcal{O}_X / \mathcal{I}(I) \to \mathcal{O}_X / 
\mathcal{I}$. Furthermore, by \cite[Theorem~3.1]{FMR} we have that
$\mathcal{O}_X / \mathcal{I}(I)$ is canonically isomorphic to
$\mathcal{O}(q^I(J_X), X/XI)$, which by Lemma~\ref{row-finite-JX} is equal to
$\mathcal{O}_{X/XI} :=
\mathcal{O}(J_{X/XI}, X/XI)$.  If we identify
$\mathcal{O}_X / \mathcal{I}(I)$ with $\mathcal{O}_{X /XI}$, then we see that
$q(\pi_{X/XI}(q^I(a))) = 0$ implies that $\pi_X(a) \in 
\mathcal{I}$ so that $a \in I$ and $q^I(a) = 0$. Thus $q$ is faithful on $\pi_{X
/ XI} (A / I)$. Furthermore, since $\mathcal{I}$ is gauge invariant, the gauge
action on $\mathcal{O}_X$ descends to an action on the quotient $\mathcal{O}_X /
\mathcal{I}$, and $q$ intertwines this action and the action on $\mathcal{O}_{X
/XI}$. Therefore Theorem~\ref{GIU} implies that $q$ is injective and
consequently $\mathcal{I}(I) = \mathcal{I}$.

To see that the above map is injective it suffices to prove that $\pi_X(a) \in 
\mathcal{I}(I)$ if and only if $a \in I$. Now $\mathcal{O}_X / \mathcal{I}(I)
$ is canonically isomorphic to $\mathcal{O}_{X / XI}$ as in the previous
paragraph. Hence $\pi_X(a) \in \mathcal{I}(I)$ implies $\pi_{A / I}( q^I(a)) =
0$, but since $\pi_{X / XI}$ is injective by Corollary~\ref{univ-repn-inj} it
follows that $q^I(a) = 0$ and $a \in I$.
\end{proof}

\begin{remark}
We mention that in \cite{MT2} we have constructed examples which show that the
above theorem does not hold if either of the hypotheses (1) or (2) are
removed.  We also mention that Katsura \cite{Kat5} has given a description of
the gauge-invariant ideals in $C^*$-algebras associated to general
$C^*$-correspondences in terms certain pairs of ideals in $A$.
\end{remark}

%%%%%%%%%%%%%%%%%%%%%%%%%%%%%%%%%%%%%%%%%%%%%%%%%%%%%%%%%%
\section{Concluding Remarks} \label{conc-sec}
%%%%%%%%%%%%%%%%%%%%%%%%%%%%%%%%%%%%%%%%%%%%%%%%%%%%%%%%%%

In Section~\ref{add-tails-sec} we gave a method for ``adding tails to
sinks" in $C^*$-correspondences; that is, given a $C^*$-correspondence $X$ we
described how to form a $C^*$-correspondence $Y$ with the property that the
left action of $Y$ is injective and $\OX$ is canonically isomorphic to a full
corner in $\OY$.  The process of adding tails to $C^*$-correspondences provides
a useful tool for extending results for augmented Cuntz-Pimsner algebras
(i.e.~$C^*$-algebras associated to $C^*$-correspondences in which $\phi$ is
injective) to $C^*$-algebras associated to general $C^*$-correspondences.  

We used this idea in Section~\ref{GIU-sec} to extend the Gauge-Invariant
Uniqueness Theorem for augmented Cuntz-Pimsner algebras to the general case.
More generally, however, we see that many questions about $C^*$-algebras
associated to correspondences may be reduced to the corresponding questions for
augmented Cuntz-Pimsner algebras.  For example, we see that for any property
that is preserved by Morita equivalence (e.g.~simplicity, AF-ness, pure
infiniteness), one need only characterize when augmented Cuntz-Pimsner algebras
will have this property, and then by adding tails one may easily deduce a
theorem for $C^*$-algebras associated to general $C^*$-correspondences.

In addition, if $p \in \mathcal{M}(\mathcal{O}_Y)$ is the projection that
determines $\mathcal{O}_X$ as a full corner of $\mathcal{O}_Y$ (so that
$\mathcal{O}_X \cong p \mathcal{O}_Y p$), then the Rieffel
correspondence from the lattice of ideals of $\mathcal{O}_Y$ to the lattice of
ideals of $\mathcal{O}_X$ takes the form $I \mapsto pIp$.  Furthermore, we
see from Lemma~\ref{proj-in-OY} that $p$ is gauge invariant, and consequently
the Rieffel correspondence preserves gauge invariance of ideals.  Thus
questions about the ideal structure of $\mathcal{O}_X$, or about
gauge-invariant ideals of $\mathcal{O}_X$, may be reduced to the corresponding
questions for ideals in the augmented Cuntz-Pimsner algebra $\mathcal{O}_Y$.

Finally, we mention that in \cite[\S4]{MT2} the method of adding tails has
proven very useful in the analysis of topological quivers.  Topological quivers,
which were first introduced in \cite[Example~5.4]{MS2}, are generalizations of
graphs in which the sets of vertices and edges are replaced by topological
spaces.  By adding tails to topological quivers in \cite{MT2} the authors are
able to reduce their analyses to the case when there are no sinks, or equivalently, to the case when the left
action of the associated $C^*$-correspondence is injective.  This simplifies
the proofs of many results for topological quivers and allows one to avoid a
number of technicalities.

\end{document}